\documentclass[11pt,leqno,twoside]{article}

\usepackage{amsfonts,amsmath,amsthm,amssymb}
\usepackage{enumerate}
\usepackage{graphics,graphicx,subfigure}

\usepackage{color}
\usepackage{todonotes}
\usepackage{cancel}
\usepackage{url}
\usepackage{hyperref}
\usepackage{makeidx}
\usepackage{showidx}
\usepackage{multicol}        
\usepackage{xspace}
\usepackage{stmaryrd}        
\usepackage{pifont}          
\usepackage{fancybox}        
\usepackage{bm}

 \usepackage{fancyhdr}
\theoremstyle{plain}
 \theoremstyle{definition}
 \newtheorem{lem}{Lemma}
 \newtheorem{defn}[lem]{Definition}
 \newtheorem{thm}[lem]{Theorem}
 \newtheorem{prop}[lem]{Proposition}
 \newtheorem{cor}[lem]{Corollary}
 \newtheorem{notn}[lem]{Notations}
 \newtheorem{pb}[lem]{Problem}
 \newtheorem{form}[lem]{Formulae}
 
 \newtheorem*{rk}{Remark}
 \newtheorem*{com}{Comment}
 \newtheorem*{ex}{Example}
 \theoremstyle{remark}

 \newcommand{\blem}{\begin{lem}}
 \newcommand{\elem}{\end{lem}}
 \newcommand{\bdefn}{\begin{defn}}
 \newcommand{\edefn}{\end{defn}}
 \newcommand{\bthm}{\begin{thm} }
 \newcommand{\ethm}{\end{thm}}
 \newcommand{\bprop}{\begin{prop}}
 \newcommand{\eprop}{\end{prop}}
 \newcommand{\bcor}{\begin{cor}}
 \newcommand{\ecor}{\end{cor}}
 \newcommand{\bnotn}{\begin{notn}}
 \newcommand{\enotn}{\end{notn}}
 \newcommand{\bpb}{\begin{pb}}
 \newcommand{\epb}{\end{pb}}
 \newcommand{\bform}{\begin{form}}
 \newcommand{\eform}{\end{form}}
 \newcommand{\brk}{\begin{rk}}
 \newcommand{\erk}{\end{rk}}
 \newcommand{\bcom}{\begin{com}}
 \newcommand{\ecom}{\end{com}}
 \newcommand{\bex}{\begin{ex}}
 \newcommand{\eex}{\end{ex}}
 \newcommand{\bpf}{\begin{proof}}
 \newcommand{\epf}{\end{proof}}

\newcommand{\cC}{\mathcal{C}}

\newcommand{\cE}{\mathcal{E}}
\newcommand{\cF}{\mathcal{F}}

\newcommand{\cI}{\mathcal{I}}

\newcommand{\cO}{\mathcal{O}}

\newcommand{\cS}{\mathcal{S}}

\newcommand{\cU}{\mathcal{U}}

\newcommand{\bC}{\mathbb{C}}

\newcommand{\bE}{\mathbb{E}}

\newcommand{\bP}{\mathbb{P}}

\newcommand{\bR}{\mathbb{R}}

\newcommand{\be}{\begin{equation}}
\newcommand{\ee}{\end{equation}}
\newcommand{\bal}{\begin{align}}
\newcommand{\eal}{\end{align}}
\newcommand{\ba}{\begin{align*}}
\newcommand{\ea}{\end{align*}}
\newcommand{\bmx}{\begin{matrix}}
\newcommand{\emx}{\end{matrix}}
\newcommand{\bbmx}{\begin{bmatrix}}
\newcommand{\ebmx}{\end{bmatrix}}
\newcommand{\bpmx}{\begin{pmatrix}}
\newcommand{\epmx}{\end{pmatrix}}
\newcommand{\bvmx}{\begin{vmatrix}}
\newcommand{\evmx}{\end{vmatrix}}

\newcommand{\wh}{\widehat}
\newcommand{\wt}{\widetilde}
\newcommand{\f}{\frac}

\newcommand{\inc}{\subseteq}

\newcommand{\setm}{\setminus}

\newcommand{\spa}{\mathrm{span}}
\newcommand{\Id}{\mathrm{Id}}

\newcommand{\tr}{\mathrm{tr}}

\newcommand{\argmin}{{\rm argmin}\,}

\newcommand{\minimize}[1]{\underset{#1}{\rm minimize}\,}

\newcommand{\la}{\lambda}
\newcommand{\La}{\Lambda}
\newcommand{\eps}{\varepsilon}
  
\newcommand{\ran}{\operatorname{range}}

\usepackage{natbib}
\title{\vspace{-20mm}
Worst-Case Completion of Tensors\\ with Approximately Few ANOVA Terms
\medskip\hrule height 1.2pt \vspace{-6mm}}
\author{Simon Foucart\footnote{S.~F. is partially funded by the NSF  (DMS-2505204).} \,
and Jingchun Shao\footnote{S.~F. and J.~S. also acknowledge support from Texas~A\&M~University through the ASCEND Initiative.} --- Texas A\&M University}
\date{\vspace{-6mm}\rule{100mm}{0.8pt}}

\newcommand\shorttitle{Completion of Tensors with Few ANOVA Terms}
\newcommand\authors{S. Foucart, J. Shao}

\begin{document}
\maketitle

\vspace{-15mm}
\begin{abstract}
In this article, 
the problem of completing a tensor from some incomplete knowledge of its entries is treated by adopting a worst-case perspective,
given the realistic assumption that the tensor's low-order ANOVA terms are dominant.
We survey and leverage some recent all-purpose results from the field of Optimal Recovery to provide solutions on a theoretical level.
But the accompanying constructions of optimal completion procedures,
which often feature semidefinite programs,
are not directly applicable in the tensor case due to the huge dimensions involved.
To resolve the issue,
we put forward a storage-friendly way to produce low-order ANOVA projections based on the fast Fourier transform (FFT),
while exploiting the specificities of the completion problem to efficiently compute regularizers and extremal eigenvalues.
Numerical experiments on synthetic tensors and real-world datasets demonstrate the accuracy and scalability of our FFT-based method.

\end{abstract}



\vspace{-5mm}
\begin{center}
\rule{100mm}{0.8pt}
\end{center}

\section{Introduction}

In this article,
we are seeking optimal algorithms for the recovery of a $d$-way tensor $X \in \bR^{n_1 \times \cdots \times n_d}$
known {\em a priori} to have approximate low complexity and 
observed {\em a posteriori} via 
the direct entries
$$
y_{\mathbf{i}} 
= X[\mathbf{i}]
= X[i_1,\ldots,i_d],
\qquad \mathbf{i} \in \cO,
$$
where $\cO$ represents the multi-index set of observed entries.
Without being too formal for now,
approximate low complexity
means that the low-order terms dominate in the ANOVA decomposition $X = \sum_{S \inc [1:d]} \Pi_S(X)$,
about which details will appear in Section~\ref{SecANOVA}.
Our recovery task consists in completing the tensor $X$ based on this knowledge and on the observations $y_{\mathbf{i}} = X[\mathbf{i}]$, $\mathbf{i} \in \cO$,
with a focus on the worst-case setting.
Thus, the task shall be carried out in the framework of Optimal Recovery.

This problem is closely related to tensor completion, where one aims at recovering a tensor from a limited number of 
entries by imposing an additional low-complexity structure. 
A common structural choice is low tensor rank, which leads to tensor analogues of matrix completion and compressive sensing. 
For example, low-rank tensor completion can be formulated as optimization programs over tensors of a fixed rank---such as Canonical Polyadic (CP) rank \citep{JainOh2014}
and Tensor Train (TT) rank 
\citep{steinlechner2016riemannian}---or 
relaxed as convex optimization programs involving the tensor nuclear norm \citep{LiuEtAl2013}. 
To solve these programs efficiently, Riemannian optimization methods have been widely developed \citep{KressnerSteinlechnerVandereycken2014}. 
Other works explore structured sampling schemes,
e.g. \citet{Zhang2016Cross} demonstrates that,
in the third-order Tucker setting,
a tensor can be exactly recovered from noiseless cross-tensor measurements whose count matches the intrinsic degrees of freedom. 
However, unlike in the matrix case, 
there is no 
universal preference among tensor norms and computations featuring these norms are notoriously difficult. 
For instance, computing the tensor rank and the tensor nuclear norm are generally NP-hard problems \citep{HillarLim2013}. 
In contrast to this traditional low-rank viewpoint, the present work assumes that the underlying tensor can be well approximated by the low-order terms of its ANOVA decomposition.

In the general framework of Optimal Recovery,
tensors $X$ give way to abstract objects $f$ living in a Banach space $F$,
the {\em a priori} information takes the form $f \in \cC$ for a so-called model set $\cC \inc F$,
and the {\em a posteriori} information occurs as $y = \La(f)$ 
for a linear map $\La: F \to \bR^m$ termed observation map.
Towards the objective of optimally recovering a quantity of interest $\Gamma(f)$ for some $\Gamma: F \to G$,
a recovery process is simply viewed as a map $\Delta: y \in \bR^m \to \Delta(y) \in G$.
To assess its performance, 
there are two important notions of worst-case errors:\vspace{-5mm}
\begin{itemize}
\item the global worst-case error of $\Delta$, as defined by
$$
{\rm gwce}(\Delta) = \sup_{f \in \cC} \|\Gamma(f) - \Delta(\La f) \|_G;
$$
\item 
the local worst-case error of $\Delta(y) =: g$
at the specific $y \in \La(\cC)$,
as defined by
$$
{\rm lwce}_y(g) = \sup_{\substack{f \in \cC \\ \La ( f )= y}} \|\Gamma(f) - g\|.\vspace{-5mm}
$$
\end{itemize}
Our goal shall be to construct optimal recovery maps $\Delta^{\rm opti}: \bR^m \to G$,
be they globally optimal,
i.e., ${\rm gcwe}(\Delta^{\rm opti}) \le {\rm gcwe}(\Delta)$ for any $\Delta: \bR^m \to G$, or locally optimal, i.e., ${\rm lcwe}_y(\Delta^{\rm opti}(y)) \le {\rm lcwe}_y(g)$ for any $g \in G$.
Of note, a locally optimal recovery map is automatically globally optimal.

A fertile specification of the general framework, 
extended to incorporate observation errors,
has been studied by the first author in the series of works \citep{foucart2022learning,foucart2023optimal,foucart2023optimalSampTA,foucart2024radius,foucart2024s,foucart2025worst}.
In this specification, 
the spaces $F$ and $G$ are Hilbert spaces,
which we assume from now on,
and the model set $\cC$ is an hyperellipsoid
$$
\cC_T := \{ f \in F: \|T(f) \| \le 1 \}
$$
relative to an operator $T: F \to F$. 
Of particular interest
is the case of an approximability set 
$$
\cC_{V,\eps} := \{ f \in F:
{\rm dist}(f,V) \le \eps\}
= \{ f \in F: \|P_{V^\perp}(f) \| \le \eps \}
$$
relative to a linear subspace $V$ of $F$ and to a parameter $\eps > 0$.
This approximability set coincides with the hyperellipsoid by taking $T = \eps^{-1} P_{V^\perp}$ to be the suitably scaled orthogonal projector  onto the orthogonal complement of $V$.
Of further particular interest is the case of an observation map $\La: F \to \bR^m$ that satisfies $\La \La^* = \Id_{\bR^m}$.

These two cases prevail in our tensor scenario.
Indeed, the low-complexity assumption translates into $X$
being well-approximated (up to some accuracy $\eps$) by its orthogonal projection onto the ANOVA space of level $s$ (for a small integer $s \ge 0$).
Moreover,
the assumption $\La \La^* = \Id_{\bR^m}$ clearly holds when $\La(X) = X_{| \cO}$ consists of entries of $X$ indexed by a set $\cO$.
Seemingly, then, 
there is not much to do,
as the above-mentioned series of works gave an (almost) complete theoretical solution of the Optimal Recovery problem at hand.
There is a caveat, though:
because of the huge nominal dimension of the space $\mathbb{R}^{n_1 \times \cdots \times n_d}$,
the computational recipes uncovered earlier
cannot be executed in a straightforward manner.
In fact, storing tensors---in a naive way---might not even be an option.
Thus, the contribution of the present work is twofold.
Firstly, it serves as a survey of disparate theoretical results---and even complements them with two new contributions.
Secondly,
and more importantly,
it puts the theoretical results to the strenuous test of numerical realization with tensors,
where the naive computational recipes must be adjusted.

The organization of the rest of the article is as follows.
In Section~\ref{SecANOVA},
in order to formally introduce our approximability set,
we recall the main features of the ANOVA decomposition,
favoring an alternative view over the traditional one.
Importantly, this view will lead to more efficient computations involving the fast Fourier transform (FFT)
and to guarantees that our worst-case completion problem is well-posed.
In Section~\ref{SecORExact}, 
we study the optimal completion problem within the scenario of accurate observations,
starting with an exposition of existing Optimal Recovery results
and continuing with their computational adaptations to the high-dimensional tensor framework.  
Section~\ref{SecORInexact} considers the optimal completion problem in the more realistic---and more difficult---scenario of inaccurate observations.
Starting again with an exposition of existing results,
it also contains new results,
notably one improving the practical construction of Chebyshev centers in a situation that comprises the tensor completion task.
It closes again with computational adaptations  to the tensor framework.
Finally, in Section~\ref{sec:numerical-experiments}, 
our numerical solutions are validated on a brief selection of relevant examples.
Their {\sc matlab} implementations can be found in the dedicated GitHub repository \url{https://github.com/Jingchun-Shao/OR_ANOVA},
which also contains a file allowing to reproduce the tests presented here.

\section{The ANOVA Decomposition}
\label{SecANOVA}

In this section,
we formalize our model assumption based on the ANOVA decomposition of a tensor.
As a matter of fact,
the ANOVA decomposition more generally applies to real-valued functions~$x$ depending on $d$ variables $t_1 \in \Omega_1, \ldots, t_d \in \Omega_d$,
each $\Omega_j$ being equipped with a probability measure~$\mu_j$.
Tensors in $\mathbb{R}^{n_1 \times \cdots \times n_d}$ are just the special case of $\Omega_j = [1 : n_j ]$ equipped with the normalized counting measure.
We will therefore now use the notation $x \in \bR^{\Omega_1 \times \cdots \times \Omega_d}$ as a compromise between generic objects $f \in F$ and proper tensors $X \in \bR^{n_1 \times \cdots \times n_d}$.

\subsection{Traditional and alternative views}

The ANOVA decomposition of a real-valued function $x \in \bR^{\Omega_1 \times \cdots \times \Omega_d}$ of $d$ variables is often encountered as the $L_2(\mu_1 \otimes \cdots \otimes \mu_d)$-orthogonal sum
\be
\label{SumANOVA}
x(t) = \sum_{S \inc [1:d]} x_S(t).
\ee 
The functions $x_S$ depend only on the variables $t_j$ for $j \in S$
and are defined recursively via
\be
\label{TradANOVA}
x_S = P_S(x) - \sum_{\substack{R \inc S\\ R \not= S}} x_R,
\ee
where the operator $P_S$ averages out all the variables outside of $S$.
Thus,
with $P_{-\ell}$ as a shorthand for $P_{[1:d] \setm \{\ell\}}$, we have 
$$
P_S = \prod_{\ell\notin S} P_{-\ell},
\qquad
\mbox{with} \quad
(P_{-\ell}x)(t_1,\dots,t_d)
:=
\int_{\Omega_\ell} x(t_1,\dots,t_{\ell-1},\tau,t_{\ell+1},\dots,t_d)\,d\mu_\ell(\tau).
$$

Besides the traditional view \eqref{SumANOVA}-\eqref{TradANOVA} of the ANOVA decomposition,
there is an alternative one found e.g. in \citep{takemura1983tensor}.
It will soon be exploited for efficient computing 
and is arguably quite simple to establish.
We include our explanation for completeness.

\bprop
The ANOVA decomposition of a $d$-variate function $x$ is
\be
\label{AltSumANOVA}
x = \sum_{S \inc [1:d]} \Pi_S(x),
\ee
where the linear operators $\Pi_S$ are mutually orthogonal projectors given by
\be
\label{AltANOVA}
\Pi_S = \prod_{j \in S} (\Id - P_{-j}) \; P_S
= \sum_{R \inc S} (-1)^{|S|-|R|} P_R.
\ee
\eprop

\bpf
The arguments only use the fact that the operators $P_{-j}
$,
which average out the $j$th variable,
are orthogonal projectors which commute with one another.
It starts with the following expansion:
$$
\Id = \prod_{i=1}^d (\Id-P_{- i} + P_{- i})
= \sum_{S \inc [1:d]} \prod_{j \in S} (\Id - P_{- j})
\prod_{\ell \not \in S} P_{- \ell}
= \sum_{S \inc [1:d]} \Pi_S.
$$
This is the decomposition \eqref{AltSumANOVA}
in which we have defined 
$$
\Pi_S
:=
\prod_{j \in S} (\Id - P_{- j})
\prod_{\ell \not \in S} P_{- \ell}
= \prod_{j \in S} (\Id - P_{-j}) \; P_S.
$$
The commutativity of the $P_{-i}$ readily implies that $\Pi_S^* = \Pi_S$ and $\Pi_S \Pi_S = \Pi_S$,
i.e., that the $\Pi_S$ are orthogonal projectors.
Moreover,
they are mutually orthogonal, 
i.e., $\Pi_S \Pi_{S'} = 0$ for $S \not= S'$---indeed, 
picking e.g. $i \in S' \setm S$, say,
one sees that $\Pi_S \Pi_{S'}$ contains the product
$P_{- i} (\Id - P_{- i}) = 0$.
Notice that the expansion of 
$\prod_{j \in S} (\Id - P_{- j})$
in the expression for $\Pi_S$ yields
$$
\Pi_S = \prod_{j \in S} (\Id - P_{- j}) \; P_S
= \sum_{R \inc S} \prod_{i \in S \setm R}(-P_{- i}) \; P_S
= \sum_{R \inc S} (-1)^{|S|-|R|} P_R,
$$ 
which is the second identity in \eqref{AltANOVA}.
Finally,
the recursive definition \eqref{TradANOVA} is retrieved by way of yet another expansion, namely  
$$
P_S = \prod_{i \in S} (\Id- P_{- i} + P_{- i}) \; P_S
= \sum_{R \inc S} \prod_{j \in R} (\Id- P_{-j}) \prod_{\ell \in S \setm R} P_{- \ell} \; P_S
= \sum_{R \inc S} \prod_{j \in R} (\Id- P_{- j}) \; P_R
= \sum_{R \inc S} \Pi_R,
$$
which yield $\displaystyle{\Pi_S = P_S - \sum_{\substack{R \inc S \\ R \not= S}} \Pi_R}$ after a rearrangement. 
\epf

\subsection{The ANOVA-based approximability set}

In order to sidestep the curse of dimensionality for the approximation of multivariate functions~$x$,
it is common to make some variable-diminution assumptions,
e.g. 
(i) $x$ belong to a weighted function space,
where the importance of each variable $t_j$ is quantified by a weight $\gamma_j > 0$ (see \citep{sloan1998quasi});
(ii) $x$ depend only on a few variables, 
either coordinate variables $t_j$ with indices $j$ belonging to an unknown  set $S$ (see \citep{devore2011approximation})
or reduced variables $\langle a_k, t \rangle$ with unknown vectors $a_k$ (see \citep{cohen2012capturing});
(iii) $x$ is partially separable, i.e., the sum of functions of few variables.
Assumption (iii) is evidently the one we shall focus on,
specifically postulating that only the first few terms of the ANOVA decomposition of $x$ matter,
say those indexed by subsets $S$ of $[1:d]$ with $|S| \le s$ for some small integer $s \ge 0$.

In many practical applications, high-dimensional functions are indeed largely governed
by individual variables and low-order interactions between them, 
i.e., somewhat equivalently, by their
low-order ANOVA terms. 
This observation motivates our use of ANOVA truncation
as a low-complexity model. 
It also helps explain the empirical success of quasi-Monte Carlo methods: as shown by \citet{Owen2003}, ANOVA
decompositions often reveal that standard quadrature test functions have an
effective dimension much smaller than the ambient dimension. 
A similar principle underlies the High-Dimensional Model Representation (HDMR) framework of 
\citep*{LiRabitz2001}, 
where complex chemical and molecular systems are approximated by sums of low-order component functions. 
From the approximation-theoretic side, \citet{Griebel2006} emphasizes that ANOVA-type
decompositions reveal the relative importance of variables and their
interactions and discusses practical examples where low-order or rapidly
decaying ANOVA structure appears, including molecular dynamics, Markov-process
representations, and high-dimensional problems in mathematical finance.

To make things more precise,
we now introduce the space 
\be
\label{DefSpaceVs}
V_s :=
\bigoplus_{|S|\le s}^\perp \ran(\Pi_S)
= \ran(\Pi_{\le s}),
\qquad \mbox{where} \quad
\Pi_{\le s} := \sum_{|S| \le s} \Pi_S.
\ee
Rather than assuming that our functions $x \in \bR^{\Omega_1 \times \cdots \times \Omega_d}$---or rather our tensors $X \in \bR^{n_1 \times \cdots \times n_d}$---belong exactly to $V_s$,
we stipulate that they are well approximated by elements of $V_s$,
say up to accuracy $\eps > 0$.
This leads us to use as model set the following approximability set:
$$
\cC_{V_s,\eps}
:= \{
X \in \bR^{n_1 \times \cdots \times n_d}:
\| X - \Pi_{\le s}(X)\| \le \eps
\},
$$
which appears as an hyperellipsoidal model set $\cC_T$ with 
$T = \eps^{-1} (\Id - \Pi_{\le s})$.
For such hyperellipsoidal model sets,
a necessary condition for the Optimal Recovery problem to even make sense is easily\footnote{If there was a nonzero $v \in \ker(T)$ such that $\La(v)=0$,
then, picking some $f_0 \in \cC_T$ with $\La(f_0)=y$,
any $f_t:= f_0 +t v$ would belong to $\cC_T$ and satisfy $\La(f_t) = 0$,
while $\|f_t - g \| \ge t \|v\| - \|f_0-g\| \to +\infty$ as $t \to + \infty$,
so the local worst-case error ${\rm lwce}_y(g)$ for the recovery of $\Gamma = \Id$ would be infinite for any $y$ and $g$.
Likewise for the global worst-case error of any $\Delta$.} seen to be $\ker(T) \cap \ker(\La) = \{0\}$.
In our situation,
this means that
\be
\label{NecCond}
V_s \cap \ker(\La) = \{0\}.
\ee
Viewed as the injectivity of $\La_{|V_s}$,
this condition implies $\dim(V_s) \le |\cO|$.
To close this subsection,
we estimate the dimension of $V_s$
before 
uncovering a specific set $\cO$ with size $|\cO| = \dim(V_s)$
for which
 Condition \eqref{NecCond} is fulfilled.

\blem
\label{LemEstimK}
The dimension of the space $V_s$ defined in~\eqref{DefSpaceVs} is
\be
\label{dimVs_exact}
k_s := \dim(V_s) = \sum_{|S| \le s} \prod_{j \in S} (n_j-1).
\ee
In the particular case $n_1 = \cdots = n_d =: n$, it can be lower- and upper-estimated as
$$
\left( \f{d \, (n-1)}{s} \right)^s
\le 
k_s \le
\left( \f{e \, d \, (n-1)}{s} \right)^s.
$$
\elem 

\bpf
Since rank and trace coincide for (not necessarily  orthogonal) projectors
and in view  of the second identity of \eqref{AltANOVA},
namely $\Pi_S = \sum_{R \inc S} (-1)^{|S|-|R|} P_R$,
we observe that, for a fixed $S \inc [1:d]$,
\begin{align*}
{\rm rk}(\Pi_S) 
& = \tr(\Pi_S)
= \sum_{R \inc S} (-1)^{|S|-|R|} \tr(P_R)
= \sum_{R \inc S} (-1)^{|S|-|R|} {\rm rk}(P_R)
= \sum_{R \inc S} (-1)^{|S|-|R|} \prod_{i\in R} n_i\\
& = \prod_{j \in S} (n_j-1).
\end{align*}
Now, since the space $V_s$ is the orthogonal sum  of the spaces $\ran(\Pi_S)$ for $|S| \le s$,
its dimension is equal to $\sum_{|S| \le s} {\rm rk}(\Pi_S)$, which justifies \eqref{dimVs_exact}.
When $n_1=\cdots=n_d=n$,
this yields
$$
\dim(V_s) = \sum_{r=0}^s \binom{d}{r} (n-1)^r.
$$
For the lower bound, we only consider the summand corresponding to $r=s$ and use $\binom{d}{r} \ge \big(\f{d}{s} \big)^s$.
For the upper~bound, we use $(n-1)^r \le (n-1)^s$ and the inequality $\sum_{r=0}^s \binom{d}{r} \le \big( \f{ed}{s} \big)^s$, see e.g. \citep[Lemma~2.6]{foucart2022mathematical}.
\epf

The exact statement for the fulfillment of Condition~\ref{NecCond} is included below,
but its proof is postponed until the next subsection,
as it uses an ingredient introduced there.

\bprop
\label{PropNecCond}
With $V_s$ being the subspace of $\bR^{n_1 \times \cdots \times n_d}$ 
defined in~\eqref{DefSpaceVs},
Condition \eqref{NecCond} holds 
for the observation map
$\La: X \in \bR^{n_1 \times \cdots \times n_d} \mapsto \big( X[\mathbf{i}], i \in \cO \big) \in \bR^{|\cO|}$  associated any multi-index set $\cO$ containing the set
$ \cI_{d,s} = 
\big\{ 
\mathbf{i}=(i_1,\ldots,i_d): |\{ j \in [1:d]: i_j \not= 1 \}| \le s 
\big\}.
$
\eprop

\subsection{Efficient computations of the ANOVA projections}

Assuming for the sake of discussion that $n_1 = \cdots = n_d = n$,
we consider a tensor $X \in \bR^{n \times \cdots \times n }$ and a subset $S$ of $[1:d]$ with size $|S| = s$.
The cost of naively computing the ANOVA projection $X_S = \Pi_S(X)$ via the recursive definition \eqref{TradANOVA},
namely via $X_S = P_S(X) - \sum_{\substack{R \inc S\\ R \not= S}} X_R$, is at least 
$2 \, s! \,  n^d$.
Indeed, evaluating $P_S(X)$ involves $n^{d-s}$ additions to average out the variables outside of $S$,
to be done at each of the $n^s$ tuples $(t_i,i \in S)$,
for a total cost of $n^d$.
Thus, with $c_r$ denoting the cost of computing $x_R$ when $|R|=r$,
we have $c_0 = n^d$,
while the recursive definition \eqref{TradANOVA} gives
$$
c_s = n^d + \sum_{r=0}^{s-1} \binom{s}{r} c_r.
$$
The claimed lower bound uses $\sum_{k \ge 0} \binom{k}{s} z^{k-s} = (1-z)^{-s-1}$ in the explicit solution of this recursion, namely
$$
c_s = n^d \,  \sum_{k=0}^\infty \f{k^s}{2^k},
\qquad \mbox{ so that } \qquad
\f{c_s}{s! \,  n^d} \ge \sum_{k=0}^\infty \binom{k}{s} \f{1}{2^k} = 2.
$$
The cost of naively computing the ANOVA projection $\Pi_{\le s}(X) = \sum_{|S| \le s} \Pi_S(X)$ will be even higher.
We are going to show that this cost can actually
be reduced to $ n^d \ln(n^d)$ through a Fourier-based approach exploiting the identity \eqref{AltANOVA} representing $\Pi_S$ as a product of orthogonal projectors.
The approach, 
conceptually identical to the cluster expansion techniques widely used in computational materials science
(see e.g. \citep{drautz_ace}),
is implemented in the dedicated repository accompanying this article.
Starting with a fixed index set $S$,
it can be formulated as follows by naturally extending the considerations to the complex setting.

\bthm
\label{ThmCheapANOVAS}
The ANOVA component of $X \in \bC^{n_1 \times \cdots \times n_d}$ indexed by $S \inc [1:d]$ can be expressed as
$$
\Pi_S(X)
= \cF^{-1} \big( M_S \odot \cF(X) \big),
$$
where $\cF$ is the $d$-dimensional Fourier transform on $\bC^{n_1 \times \cdots \times n_d }$ producing discrete Fourier coefficients 
and where the mask $M_S \in \{0,1\}^{n_1 \times \cdots \times n_d}$
is given for $\mathbf{i} = (i_1,\ldots,i_d) \in [1:n_1] \times \cdots \times [1,n_d]$ by 
$$
(M_S)_{\mathbf{i}}
= \begin{cases}
1 & \mbox{ if }
i_j \not= 1 \mbox{ for } j \in S  \mbox{ and } i_\ell = 1 \mbox{ for } \ell \not\in S,\\
0 & \mbox{ otherwise}.
\end{cases}
$$
\ethm

\bpf
For each $j \in [1:d]$,
let $(\varphi^{(j)}_1, \varphi^{(j)}_2, \ldots, \varphi^{(j)}_{n_j})$ denote an unnormalized Fourier basis of $\bC^{n_j}$
with $\varphi^{(j)}_{1}$ being equal to the all-one vector.
Then,
for $\mathbf{i} = (i_1,\ldots,i_d) \in [1:n_1] \times \cdots \times [1,n_d]$,
let $\phi_{\mathbf{i}} := \varphi^{(1)}_{i_1} \otimes \cdots \otimes \varphi^{(d)}_{i_d}$.
It is easily seen---recalling that $P_{-\ell}$ averages out $\ell$-th variable---that
$$
P_{-\ell}(\phi_{\mathbf{i}})
= \begin{cases}
\phi_{\mathbf{i}} & \mbox{ if } i_\ell = 1,\\
0 & \mbox{ if } i_\ell \not= 1.
\end{cases}
$$
Therefore,
in view of the leftmost identity of \eqref{AltANOVA} affirming that $\Pi_S = \prod_{j \in S} (\Id - P_{-j}) \; \prod_{\ell \not\in S} P_{- \ell}$,
we obtain
$$
\Pi_S(\phi_{\mathbf{i}})
= \left.
\begin{cases}
\phi_{\mathbf{i}} & \mbox{ if } i_j \not= 1 \mbox{ for all }
j \in S \mbox{ and } i_\ell = 1 \mbox{ for all } \ell \not\in S,  \\
0 & \mbox{ otherwise}.
\end{cases}
\right\} = (M_S)_{\mathbf{i}} \, \phi_{\mathbf{i}}.
$$
Finally, for any $X \in \bC^{n_1 \times \cdots \times n_d}$ written as $X = \cF^{-1}\big( \cF(X) \big)
= \sum_{\mathbf{i}} \cF(X)_{\mathbf{i}} \, \phi_{\mathbf{i}}$,
we conclude that
$$
\Pi_S(X) = 
\sum_{\mathbf{i}} \cF(X)_{\mathbf{i}} \, (M_S)_{\mathbf{i}} \, \phi_{\mathbf{i}}
= \sum_{\mathbf{i}}  (M_S \odot \cF(X))_{\mathbf{i}} \, \phi_{\mathbf{i}}
= \cF^{-1} \big( M_S \odot \cF(X) \big),
$$
as announced.
\epf

Our claim about the cost of the above approach stems from neglecting the cost of constructing~$M_S$,
which is precomputed offline,
and considering only the costs of constructing the $d$-dimensional Fourier transform on $\bC^{n^d}$ and its inverse,
each being known (see e.g. \citep[p.76]{dudgeon1984multidimensional}) to be of order $n^d \ln(n^d) $.
In addition,
recall that we are more interested in computing the orthogonal projection of $X$ onto the approximation space $V_s$,
i.e., in $\Pi_{\le s}(X) = \sum_{|S| \le s} \Pi_S(X)$,
than in computing an individual $\Pi_S(X)$.
The computational cost remains of order $n^d \ln(n^d) $,
though.
This is because 
the following corollary indicates that the majority of the cost still stems from the Fourier transform and its inverse.

\bcor
\label{CorCompPi<=s}
The order-$s$ ANOVA projection of $X \in \bC^{n_1 \times \cdots \times n_d}$ can be expressed as
\be
\label{CheapAction}
\Pi_{\le s}(X)
= \cF^{-1} \big( M_{\le s} \odot \cF(X) \big),
\ee
where the binary mask $M_{\le s} \in \{0,1\}^{n_1 \times \cdots \times n_d}$ is precomputed offline
through its entries given for $\mathbf{i} = (i_1,\ldots,i_d) \in [1:n_1] \times \cdots \times [1:n_d]$ by
$$
(M_{\le s})_{\mathbf{i}}
= \begin{cases}
1 & \mbox{ if }
|\{ j \in [1:d]: i_j \not= 1 \}| \le s,\\
0 & \mbox{ otherwise}.
\end{cases}
$$
\ecor

\bpf
In view of $\Pi_{\le s}(X)= \sum_{|S| \le s} \Pi_S(X)$,
the identity~\eqref{CheapAction} follows with $M_{\le s} := \sum_{|S| \le s} M_S$.
Note that $M_{\le s}$ takes values in $\{0,1\}$
because the summands $M_S \in \{0,1\}^{n_1 \times \cdots \times n_d}$ are disjointly supported,
as Theorem~\ref{ThmCheapANOVAS}
says that $(M_{S})_{\mathbf{i}} = 1$ if and only if
$S = \{ j \in [1:d]: i_j \not= 1 \}$.
Moreover,
the entry $(M_{\le s})_{\mathbf{i}}$ equals $1$
if and only if there exists $S \inc [1:d]$ with $|S| \le s$ such that $(M_S)_{\mathbf{i}} = 1$,
i.e., such that $S = \{ j \in [1:d]: i_j \not= 1 \}$.
In other words, one has $(M_{\le s})_{\mathbf{i}} = 1$
if and only if $|\{ j \in [1:d]: i_j \not= 1 \}| \le s$,
as desired.
\epf 

As a theoretical consequence of the above interpretation of the order-$s$ ANOVA projection,
we now prove that Condition~\ref{NecCond}
can be fulfilled with an observation set of minimal size $k_s = \dim(V_s)$.
This set is the support of $M_{\le s}$,
denoted by $\cI_{d,s}$ in Proposition~\ref{PropNecCond}.
However, for notational convenience,
the proof below shall 
feature a version of the Fourier transform where the all-one vector appear in the last position.
As a result, the necessary cyclic shifts transform this set into 
$$
\cI_{d,s} = 
\big\{ 
\mathbf{i}=(i_1,\ldots,i_d): |\{ j \in [1:d]: i_j \not= n_j \}| \le s 
\big\}.
$$

\bpf[Proof of Proposition~\ref{PropNecCond}]
Showing that  $V_s \cap \ker(\La) = \{ 0 \}$
consists in proving that,
if $X \in {\rm range}(\Pi_{\le s})$ satisfies $X[\mathbf{i}] = 0$ for all $\mathbf{i} \in \cO$,
then $X = 0$.
In view of~\eqref{CheapAction},
such an $X$ can be written as $X = \cF^{-1}(Z)$,
where $Z := M_{\le s} \odot \cF(X)$ is supported on $\cI_{d,s}$.
Taking $\cO \supseteq \cI_{d,s}$ into account,
it is therefore enough to prove that
\be
\label{IndHyp}
\mbox{if } Z \mbox{ is supported on } \cI_{d,s}
\mbox{ and satisfies }
\cF(Z)[\mathbf{i}] = 0
\mbox{ for all } i \in \cI_{d,s},
\mbox{ then }
Z = 0.
\ee
We shall prove by induction on $d \ge 1$ that \eqref{IndHyp} holds for all $s=1,\ldots,d$.
The base case $d=1$ is nothing but the invertibility of the one-dimensional discrete Fourier transform.
Now, assuming that \eqref{IndHyp} holds up to $d-1$ for all relevant~$s$,
our goal is to prove that it holds for $d$ and  all relevant~$s$.
As the case $s=d$ is again nothing but the invertibility of the $d$-dimensional discrete Fourier transform, we suppose that $s \le d-1$.
At this point, it is useful to remark that
$$
\cI_{d,s} = (\cI_{d-1,s} \times \{n_d\}) 
\sqcup (\cI_{d-1,s-1} \times [1:n_d-1])
\qquad \mbox{and} \qquad
\cI_{d-1,s-1} \inc \cI_{d-1,s}.
$$
For instance,
the former is used to explicit the identities $0 = \cF(Z)[\mathbf{i}]$ for any $\mathbf{i} \in \cI_{d,s}$ as
\begin{align*}
0 & = 
\sum_{\mathbf{k} \in \cI_{d,s}} 
Z[k_1,\ldots,k_{d-1},k_d] 
\; 
e^{\mathrm{i} 2 \pi \left( \f{i_1 k_1}{n_1} + \cdots + \f{i_{d-1}k_{d-1}}{n_{d-1}} + \f{i_{d} k_{d}}{n_{d}} \right) } \\
& = \sum_{ (k_1,\ldots,k_{d-1}) \in \cI_{d-1,s}} 
Z[k_1,\ldots,k_{d-1},n_d] 
\; 
e^{\mathrm{i} 2 \pi \left( \f{i_1 k_1}{n_1} + \cdots + \f{i_{d-1}k_{d-1}}{n_{d-1}} \right) }\\
& + 
\sum_{(k_1,\ldots,k_{d-1}) \in \cI_{d-1,s-1}}
\bigg(
\sum_{k_d = 1}^{n_d - 1} Z[k_1,\ldots,k_{d-1},k_d] 
e^{\mathrm{i} 2 \pi \f{i_d k_d}{n_d}}
\bigg)
\;
e^{\mathrm{i} 2 \pi \left( \f{i_1 k_1}{n_1} + \cdots + \f{i_{d-1}k_{d-1}}{n_{d-1}} \right) } .
\end{align*}
At this point, there are two cases to separate,
namely $\mathbf{i} \in  \cI_{d-1,s} \times \{n_d\}$ and $\mathbf{i} \in \cI_{d-1,s-1} \times [1:n_d-1]$.
In the fist case, one obtains that,
for all $(i_1,\ldots,i_{d-1}) \in \cI_{d-1,s}$,
\begin{align*}
0 & = 
\sum_{ (k_1,\ldots,k_{d-1}) \in \cI_{d-1,s} \setm \cI_{d-1,s-1}} 
Z[k_1,\ldots,k_{d-1},n_d] 
\; 
e^{\mathrm{i} 2 \pi \left( \f{i_1 k_1}{n_1} + \cdots + \f{i_{d-1}k_{d-1}}{n_{d-1}} \right) }\\
& + 
\sum_{(k_1,\ldots,k_{d-1}) \in \cI_{d-1,s-1}}
\bigg(
\sum_{k_d = 1}^{n_d} Z[k_1,\ldots,k_{d-1},k_d] 
\bigg)
\;
e^{\mathrm{i} 2 \pi \left( \f{i_1 k_1}{n_1} + \cdots + \f{i_{d-1}k_{d-1}}{n_{d-1}} \right) } .
\end{align*}
Invoking the induction hypothesis 
for $(d-1,s)$,
one deduces that
\begin{align}
\label{PfNecCond1}
& \mbox{for all } (k_1,\ldots,k_{d-1}) \in \cI_{d-1,s} \setm \cI_{d-1,s-1},
& Z[k_1,\ldots,k_{d-1},n_d] & = 0,\\
\label{PfNecCond2}
& \mbox{for all } (k_1,\ldots,k_{d-1}) \in \cI_{d-1,s-1},
&  \sum_{k_d = 1}^{n_d} Z[k_1,\ldots,k_{d-1},k_d] 
& = 0.
\end{align}
In the second case, 
making use of \eqref{PfNecCond1},
one obtains that,
for all $(i_1,\ldots,i_{d-1}) \in \cI_{d-1,s-1}$ and all $i_d \in [1:n_d - 1]$,
$$
0 = 
\sum_{(k_1,\ldots,k_{d-1}) \in \cI_{d-1,s-1}}
\bigg(
\sum_{k_d = 1}^{n_d } Z[k_1,\ldots,k_{d-1},k_d] 
e^{\mathrm{i} 2 \pi \f{i_d k_d}{n_d}}
\bigg)
\;
e^{\mathrm{i} 2 \pi \left( \f{i_1 k_1}{n_1} + \cdots + \f{i_{d-1}k_{d-1}}{n_{d-1}} \right) } .
$$
Invoking the induction hypothesis 
for $(d-1,s-1)$,
one deduces that, for all $(k_1,\ldots,k_{d-1}) \in~\cI_{d-1,s-1}$,
$$
\sum_{k_d = 1}^{n_d } Z[k_1,\ldots,k_{d-1},k_d] 
e^{\mathrm{i} 2 \pi \f{i_d k_d}{n_d}}
= 0
\qquad \mbox{for all } i_d \in [1:n_d - 1],
$$
which also holds for $i_d = n_d$ by virtue of~\eqref{PfNecCond2}.
Collectively interpreted as the one-dimensional discrete Fourier transform of $Z[k_1,\ldots,k_{d-1},\cdot]$,
these equalities imply that 
$Z[k_1,\ldots,k_{d-1},k_d] = 0$
for all $(k_1,\ldots,k_{d-1}) \in \cI_{d-1,s-1}$ and all $k_d \in [1:n_d]$.
In particular,
the latter yields $Z[\mathbf{k}] = 0$ for all $\mathbf{k} \in \cI_{d-1,s-1} \times [1:n_d-1]$.
It also yields $Z[k_1,\ldots,k_{d-1},n_d] = 0$ for all $(k_1,\ldots,k_{d-1}) \in \cI_{d-1,s-1}$,
which,
combined with~\eqref{PfNecCond1}, implies that $Z[k_1,\ldots,k_{d-1},n_d] = 0$ for all $(k_1,\ldots,k_{d-1}) \in \cI_{d-1,s}$,
or equivalently that $Z[\mathbf{k}] = 0$ for all $\mathbf{k} \in \cI_{d-1,s} \times \{n_d\}$.
Altogether, these two cases ensures that $Z[\mathbf{k}] = 0$ for all $\mathbf{k} \in \cI_{d,s}$.
Since $Z$ was supported on $\cI_{d,s}$,
one concludes that $Z = 0$,
hence establishing \eqref{IndHyp} for $(d,s)$.
This concludes the inductive proof.
\epf

\section{Optimal Completion from Accurate Observations}
\label{SecORExact}

In this section,
we discuss the Optimal Recovery problem,
as abstractly described in the introduction
but specified to the Hilbert setting and the hyperellipsoidal model set $\cC_T = \{ f \in F : \|T(f)\| \le 1 \}$.
In particular, we continue---for just a little longer---to work under the idealized assumption that the vector $y = \La(f)$ is not corrupted by any observation error.

\subsection{Summary of known results}

\paragraph{Fixed worst-case errors.} 

For a fixed and linear recovery map $\Delta^{\rm lin}: \bR^m \to G$,
determining its {\em global} worst-case error is rather simple.
Its square can be expressed as
\begin{align*}
{\rm gwce}(\Delta^{\rm lin})^2
& := \sup_{\|T(f) \|^2 \le 1} 
\|\Gamma(f) - \Delta^{\rm lin}(\La f)\|^2\\
& = \inf_{c \ge 0} \; c  \quad \mbox{s.to }
\|(\Gamma - \Delta^{\rm lin} \La)(f)\|^2 \le c \, \|T(f)\|^2
\; \; \forall f \in F.
\end{align*}
When the Hilbert space $F$ is of finite dimension, say  $N$, the above translates into an expression computable by semidefinite programming\footnote{If $T$ was invertible (which is not the case for an approximability set), then solving a semidefinite program is an overkill: an eigenvalue computation is enough, as ${\rm gwce}(\Delta^{\rm lin}) = {\rm eig}_{\max}((\Gamma - \Delta^{\rm lin} \La) T^{-1})$.}, namely
\be
\label{GWCE^2}
{\rm gwce}(\Delta^{\rm lin})^2
= \inf_{c \ge 0} \; c \quad \mbox{s.to }
c \, T^*T - (\Gamma - \Delta^{\rm lin} \La)^*(\Gamma - \Delta^{\rm lin} \La) \succeq 0.
\ee
In the more delicate {\em local} setting,
for fixed $y \in \bR^m$ and $g \in G$,
with 
$$
f_y := \argmin \{ \|T (f)\|: \; \La (f) = y \},
$$
it was shown in Theorem 1 of \citep*{foucart2022learning}
(albeit formulated slightly differently and only for $\Gamma = \Id$)
that the squared local worst-case error can be expressed as
\begin{align*}
{\rm lwce}_y(g)^2
& : = \sup_{\substack{\|T (f)\|^2 \le 1 \\ \La (f) = y}} \|\Gamma(f) - g\|^2\\
&  = \inf_{c \ge 0, d \in \bR}
\|\Gamma(f_y) - g\|^2 + d\\
& \phantom{==} \mbox{ s.to }
\langle ( cT^*T - \Gamma^* \Gamma ) h,h \rangle 
+ 2 \langle \Gamma^*(g - \Gamma(f_y)), h \rangle
+ d  - c(1 - \| T(f_y) \|^2) \ge 0
\; \forall h \in \ker(\La).
\end{align*}
The main tool for the derivation of this expression was the S-lemma.
When the Hilbert space $F$ is of finite dimension~$N$,
considering a matrix $M \in \bR^{N \times (N-m)}$ whose columns span $\ker(\La)$,
the above again translates into an expression computable by semidefinite programming, namely 
\be
\label{LWCE^2}
{\rm lwce}_y(g)^2
= 
\|\Gamma(f_y) - g\|^2
+ \inf_{c \ge 0, d \in \bR}
d: \quad
\bbmx
M^*(c T^* T - \Gamma^* \Gamma  ) M & \vline & M^* \Gamma^* (g - \Gamma f_y )\\
\hline 
(g - \Gamma f_y)^* \Gamma M & \vline & d - c(1  - \|T f_y \|^2)
\ebmx \succeq 0.
\ee

\paragraph{Minimized worst-case errors.} 
In both the global and local settings,
computing an optimal recovery map becomes possible by  further minimization.
In the global setting, admitting the theoretical fact that there are linear recovery maps among the optimal ones,
we can minimize the expression \eqref{GWCE^2} over all linear maps $\Delta^{\rm lin}$.
We indeed end up with a computationally feasible program after noticing that the constraint $c \, T^*T - (\Gamma - \Delta^{\rm lin} \La)^*(\Gamma - \Delta^{\rm lin} \La) \succeq 0$ can be reformulated as
$$
\bbmx
\Id & \vline & \Gamma - \Delta^{\rm lin} \La\\
\hline 
(\Gamma - \Delta^{\rm lin} \La)^* & \vline & c \, T^* T 
\ebmx \succeq 0.
$$
In the local setting,
for each $y \in \La(\cC_T)$,
finding the minimizer\footnote{This minimizer is known as a Chebyshev center (of $\{f \in F: \|T(f)\|\le 1, \La(f)=y\}$) and is unique in our~situation.}  over all $g$ of  the expression \eqref{LWCE^2}  is again computationally feasible,
thanks to the affine dependence on $(c, d, g)$ of the constraint and to the remark that 
$$
\|\Gamma(f_y) - g\|^2 = \inf_{t \ge 0} 
\; t : \quad \bbmx
\Id & \vline & \Gamma(f_y) - g\\
\hline
(\Gamma(f_y) - g)^* & \vline & t
\ebmx \succeq 0,
$$
where the latter constraint depends affinely on $(g,t)$.
All in all, one arrives at a semidefinite program with variables $c$, $d$, $g$, and $t$.

We shall not discuss efficient ways to solve such semidefinite programs.
In fact, we would rather bypass them,
since they are quite inappropriate to deal with tensors. 
Indeed, current semidefinite solvers struggle when the involved matrices have size in the thousands---not that they take a long time to run, but simply that they do not run because of memory issues.
To the point, in the present situation, we are faced  with matrices of size equal to the dimension of the tensor space,
i.e., $n_1 \times \cdots \times n_d$,
so already ten thousands for four-way tensors with each mode dimension equal to ten.
Thus, an alternative to this semidefinite programming approach would be welcome,
and fortuitously there is one.
It covers both the global and local settings
because one explicitly constructs the locally optimal recovery map outputting the Chebyshev center---usually called central algorithm in Information-Based Complexity---which is then automatically globally optimal.
For $\Gamma  = \Id$,
it is known (see \citep{binev2017data} for the case $T = \eps^{-1} P_{V^\perp}$ and even \citep[Theorem~299.1]{kowalski1995selected} for an arbitrary $T$)
to be given via the minimal-norm interpolant (aka spline algorithm) as
\be
\label{CheCenterViaOpt}
\Delta^{\rm cheb}:
y \in \bR^m 
\mapsto
f_y = \Big[ \underset{f \in F}{\argmin} \|T(f)\|:
\;  \La(f)= y \Big] \in F.
\ee
This turns out to be a linear map, with expression $\Delta^{\rm cheb} = (T^* T)^{-1} \La^* \big( \La (T^* T)^{-1} \La^* \big)^{-1}$ when $T$ is invertible.
For an arbitrary linear $\Gamma$, 
it can be said 
(see Lemma 6 in \citep{foucart2025worst})
that $\Gamma \circ \Delta^{\rm cheb}$ is locally, hence globally, optimal and that the minimal squared local worst-case error equals
\be
\label{CheRadEigMax}
\min_{g}{\rm lwce}_y(g)^2
= {\rm eig}_{\max}\big( [T_{|W}^* T_{|W}]^{-1/2} \Gamma_{|W}^* \Gamma_{| W} [T_{|W}^* T_{|W}]^{-1/2} \big) \times 
\big( 1 - \|T (f_y)\|^2 \big),
\ee
where ${\rm eig}_{\max}$ denotes the largest eigenvalue and where $W$ stands for $\ker(\La)$.
Note that the condition $\ker(T) \cap \ker(\La) = \{0\}$ is implicitly used here.

We will discuss below efficient ways to compute \eqref{CheCenterViaOpt} and \eqref{CheRadEigMax} for our tensorial problem---certainly, the explicit expression of $\Delta^{\rm cheb}$ is not suitable, 
since $T^* T$ is too large for direct inversion, 
having size $N \times N$ with $N = n_1  \cdots  n_d$.
But before specifying to tensors,
we point out another expression,
which is valid in a general setting.
It involves the Riesz representers $(u_1,\ldots,u_m)$ of the coordinate functionals of $\La: F \to \bR^m$,
so that $\La(f) = [\langle u_1,f \rangle; \ldots; \langle u_m,f \rangle]$ for all $f \in F$,
as well as a basis $(v_1,\ldots,v_k)$ of $V$.
These appear though the  Gramian $G \in \bR^{m \times m}$
of $(u_1,\ldots,u_m)$
and the cross Gramian $C \in \bR^{m \times k}$ 
of $(u_1,\ldots,u_m)$ with $(v_1,\ldots,v_k)$,
which have entries
$$
G_{i,i'} = \langle u_i, u_{i'} \rangle ,
\qquad
C_{i,j} = \langle u_i, v_j \rangle ,
\qquad i,i'=1,\ldots,m,
\; j = 1,\ldots,k.
$$
The desired expression was obtained in \citep*[Theorem~2]{foucart2022learning}
(see also \citep[Proposition 10.3]{foucart2022mathematical})
and reads
\be
\label{FlaGram}
\Delta^{\rm cheb}(y)
= \sum_{i=1}^m a_i u_i + \sum_{j=1}^k b_j v_j,
\ee
where 
$$
b = (C^\top G^{-1} C)^{-1} C^\top G^{-1} y
\qquad \mbox{and} \qquad
a = G^{-1}( y -  Cb ).
$$
Notice that the sizes of the matrices to invert have gone down from the dimension $N$ of the ambient space to the more manageable dimensions $m$ of $\ran(\La)$ and $k$ of $V$.
Furthermore, in our tensor completion problem---and any problem where $\La \La^* = \Id$---the system $(u_1,\ldots,u_m)$ is an orthonormal basis of $\ran(\La^*)$, so
no inversion of $G$ is required since then $G = I_m$.

\subsection{Tensor computations}

At first sight, the summary we just presented seems to solve our abstract Optimal Recovery problem in the accurate scenario,
at least in theory.
In practice,
constructing the optimal recovery map of~\eqref{CheCenterViaOpt} and evaluating the minimal worst-case error of~\eqref{CheRadEigMax} can run into computational challenges,
especially when the dimension $N$ of the ambient space is large.
This is definitely the case for our tensor completion problem,
where the dimension $N = n_1 \cdots n_d$ grows exponentially with $d$.
Another issue is storage:
when $N$ is that large,
storing the operator $T$ in matrix form is unsuitable
and we should instead exploit affordable ways to apply $T$.
Thankfully,
in case $T = \eps^{-1} P_{V_s^\perp}$, which is a multiple of  $\Id - \Pi_{\le s}$, the recipe \eqref{CheapAction} provides just that.
The discussion below focuses on this~case.

\paragraph{Optimal recovery map.}
As already mentioned,
even if $T$ was invertible,
it is unrealistic to use direct methods for matrix inversions in the construction of the optimal recovery map via the formula $\Delta^{\rm cheb} = (T^* T)^{-1} \La^* \big( \La (T^* T)^{-1} \La^* \big)^{-1}$.
Indirect methods are more appropriate,
particularly the conjugate gradient method,
which is the method of choice when the operator is positive semidefinite and accessible only through its action.
The method essentially aims at solving the minimization problem~\eqref{CheCenterViaOpt} using iterative solvers.
As these solvers handle unconstrained problems more naturally than constrained one,
we approximate the exact solution by way of a regularized problem that adds a large penalty enforcing the constraint to the objective.
More details on this regularization will be given in Subsection~\ref{SubSecTensorCompInac},
which will reveal that,
provided $\La \La^* = \Id$,
the exact---not approximate---solution to \eqref{CheCenterViaOpt} can be obtained by solving two regularized problems.

Our preferred strategy is to exploit the formula~\eqref{FlaGram}.
In our tensor completion problem,
the basis $(u_1,\ldots,u_m)$ of ${\rm range}(\La^*)$ is given by the system $(E_{\mathbf{i}}, \mathbf{i} \in \cO)$,
where $E_{\mathbf{i}} \in \bR^{n_1 \times \cdots \times n_d}$ is the canonical tensor equal to one on the $\mathbf{i}$th entry and to zero on all other entries,
and so $G = \big[ \langle u_i,u_{i'} \rangle \big] = I_m$.
As for the basis $(v_1,\ldots,v_{k_s})$ of $V_s$,
Corollary~\ref{CorCompPi<=s}
implies that it can be taken as the orthonormal system 
$(\cF^{-1}(E_{\mathbf{j}}), \mathbf{j} \in \cI_{d,s})$.
The computational cost of producing this system is $k_s \, n^d \ln(n^d)$.
Ignoring the cost of forming $\sum_i a_i u_i + \sum_j b_j v_j$ when $a$ and $b$ are available,
the remaining cost for realizing the formula~\eqref{FlaGram} then divides as:\vspace{-5mm}
\begin{enumerate}

\item[(i)] computing the cross Gramian $C \in \bR^{m \times k_s}$:
this is negligible, because the $\cF^{-1}(E_{\mathbf{j}})$ have been produced above, so that  $C = \big[  \langle u_i, v_j \rangle \big]= \big[ \cF^{-1}(E_{\mathbf{j}})[\mathbf{i}]  \big]$ is already accessible;

\item[(ii)] computing the vector $b \in \bR^{k_s}$:
forming $C^\top y$ requires $\cO(k_s m)$ operations,
forming $C^\top C$ requires $\cO(k_s^2 m)$,
and computing $b = (C^\top C)^{-1} C^\top y$ by solving the linear system $(C^\top C) b = C^\top y$ requires $\cO(k_s^3)$,
contributing to a total number of
$\cO(k_s m + k_s^2 m + k_s^3) = \cO(k_s^2 m)$ operations,
in view of $k_s \le m$;

\item[(iii)] computing the vector $a \in \bR^m$:
forming $a = y - Cb$ only requires a total number of
$\cO(k_s m)$ operations, which is dominated by the cost of (ii).
\vspace{-5mm}
\end{enumerate}

Putting everything together,
the total cost of computing $\Delta^{\rm cheb}(y)$ is
$$
\cO( k_s (n^d \ln(n^d) + k_s m) ) \quad \mbox{ operations}.
$$
For instance, 
thinking of $s$ as a fixed  small integer,
so that $k_s \asymp d^s n^s$ by Lemma~\ref{LemEstimK},
and in the realistic situation where $m \asymp k_s$,
we have $k_s m \lesssim k_s^2 \lesssim d^{2s} n^{2s} \lesssim n^d$,
resulting in a cost of $\cO(k_s n^d \ln(n^d))$. 
This is not much larger than the cost $\cO(n^d \ln(n^d) )$ of computing one ANOVA projection.

In passing, we mention another strategy for the computation of $\Delta^{\rm cheb}(y) = f_y$.
The spirit is that of the representer theorem:
since the solution has the form \eqref{FlaGram},
we can exploit this knowledge when performing the minimization of \eqref{CheCenterViaOpt},
leading to the optimization program
$$
\underset{a \in \bR^m, b \in \bR^k}{\minimize \; }
\bigg\| \sum_{i=1}^m a_i P_{V^\perp} u_i \bigg\|
\qquad \mbox{s.to } \;
a_i + \sum_{j=1}^k b_j \langle u_i,v_j \rangle = y_i 
\quad \mbox{for all }  i \in [1:m].
$$
We did not explore if there is a computational advantage in solving the latter via iterative methods,
say via conjugate gradient.

\paragraph{Minimal worst-case error.}
Once we have computed the Chebyshev center $f_y$---that is, the locally optimal recovery map evaluated at $y$---determining the Chebyshev radius---that is, the minimal local worst-case error of $f_y$ at $y$---is not necessarily straightforward.
Resorting to \eqref{LWCE^2} with $g=f_y$ is natural, of course,
but this involves a semidefinite program whose matrices are too large to handle in a tensor setting.
Resorting to \eqref{CheRadEigMax} is more promising:
as $f_y$ has just been obtained,
computing the term $(1 - \|T(f_y)\|^2)$ is not a problem,
so it remains to compute the largest eigenvalue of the positive semidefinite operator $B^{-1/2} A B^{-1/2}$,
where $A:=\Gamma_{|W}^* \Gamma_{| W}$ and $B:= T_{|W}^* T_{|W}$.
Again, it is unrealistic to store these operators in matrix form,
but their actions can be computed relatively cheaply,
say when $T = \eps^{-2} P_{V_s^\perp} = \eps^{-2}(\Id - \Pi_{\le s})$.
In this case, the power method is appropriate.
It consists (usually with $B=I$) of iterating the scheme $B x^{(t)} = A x^{(t-1)}$, so that $\langle x^{(t)}, x^{(0)} \rangle  / \langle x^{(0)}, x^{(0)} \rangle \to {\rm eig}_{\max}(B^{-1} A ) = {\rm eig}_{\max} (B^{-1/2} A B^{-1/2})$ as $t \to \infty$.
A more in-depth discussion will appear in Subsection~\ref{SubSecTensorCompInac}.

\section{Optimal Completion from Inaccurate Observations}
\label{SecORInexact}

In this section,
we discuss the Optimal Recovery problem in the more realistic situation of an imperfect observation process. 
In an abstract formulation,
the {\em a posteriori} information is now given as $y = \La(f) + e$ for some unknown error vector $e \in \bR^m$.
As for the {\em a priori} information,
it still reads $f \in \cC$ for some model set $\cC \inc F$,
but is also augmented with $e \in \cE$ for some uncertainty set $\cE \inc \bR^m$.
For the recovery of a quantity of interest $\Gamma: F \to G$,
the two notions of worst-case error are adjusted to 
\begin{align*}
{\rm gwce}(\Delta)
& = \sup_{\substack{f \in C\\ e \in \cE}} \|\Gamma(f) - \Delta(\La f + e)\|,\\
{\rm lwce}_y(g)
& = \sup_{\substack{ f \in \cC\\ y - \La f \in \cE}} \|\Gamma(f) - g \|.
\end{align*}
We note that this inaccurate scenario can be reduced to the accurate scenario by considering the compound object $\wt{f} := (f,e)$,
with {\em a priori information} $\wt{f} \in \wt{\cC}:= \cC \times \cE$
and {\em a posteriori information} $y = \wt{\La}(\wt{f}) := \La f + e$.
However,
this formal reduction is of limited use for the construction of optimal recovery maps.
As a case in point,
assuming both model set $\cC$ and uncertainty set $\cE$ to be hyperellipsoids,
the compound model set $\wt{\cC}$ is not an hyperellipsoid anymore,
so the results from the previous section are inapplicable.
Below, we (almost) untangle this case,
with specific model and uncertainty sets given as
\begin{align*}
\cC_T 
& = \{ f \in F: \, \|T(f) \| \le 1 \},  \\
\cE_U 
& = \{ e \in \ell_2^m: \, \|U(e)\| \le 1 \},
\end{align*}
for some operators $T: F \to F$ and $U: \ell_2^m \to \ell_2^m$.
Of particular interest are $T = \eps^{-1} P_{V^\perp}$ and $U = \eta^{-1} \Id$,
so that the {\em a priori } information reads ${\rm dist}(f,V) \le \eps$ and $\|e\| \le \eta$.

\subsection{Summary of known results}

\paragraph{Fixed worst-case errors.}
For a fixed {\rm linear} recovery map $\Delta^{\rm lin}: \bR^m \to G$,
the squared {\em global} worst-case error is
\begin{align}
\nonumber
{\rm gwce}(\Delta^{\rm lin})^2
& := \sup_{\substack{\|T f \|^2 \le 1 \\ \|U e\|^2 \le 1}} \| \Gamma(f) - \Delta^{\rm lin} (\La f + e) \|^2\\
\label{GWCE-lin}
& = \inf_{a,b \ge 0} a+b
\quad \mbox{s.to }
a \|Tf\|^2 + b  \|U e\|^2 \ge
\|(\Gamma - \Delta^{\rm lin} \La)f - \Delta^{\rm lin} e\|^2
\quad \forall f \in F, e \in \bR^m.
\end{align}
Although this has not been stated explicitly before
(except in \citep[Lemma~12]{foucart2023optimal} for $\Gamma = \Id$, $T = \eps^{-1} P_{V^\perp}$, and $U = \eta^{-1} \Id$),
the argument simply relies on (the infinite-dimensional version of) Polyak S-procedure involving three quadratic forms. 
When the Hilbert space $F$ is of finite dimension, the above translates into an expression computable by semidefinite programming, namely
$$
{\rm gwce}(\Delta^{\rm lin})^2
= \inf_{a,b \ge 0} a+b
\; \mbox{s.to }
\bbmx
a T^* T & \vline & 0\\
\hline 
0 & \vline & b U^* U
\ebmx
\succeq 
\bbmx
(\Gamma - \Delta^{\rm lin} \La)^* (\Gamma - \Delta^{\rm lin} \La) & \vline & (\Gamma - \Delta^{\rm lin} \La)^* \Delta^{\rm lin}\\
\hline
(\Delta^{\rm lin})^* (\Gamma - \Delta^{\rm lin} \La)
& \vline & (\Delta^{\rm lin})^* \Delta^{\rm lin}.
\ebmx .
$$
In the more delicate {\em local} setting,
because of the presence of linear terms in the quadratic functionals 
(e.g. $\|\Gamma(f) - g\|^2$ contains $\langle \Gamma(f),g \rangle$),
the S-procedure is not exact,
so a correct expression for the worst-case error is not immediately available---only an upper bound is.

\paragraph{Minimized global worst-case error.}
In the less challenging {\em global} setting,
there is now a complete solution for the construction of a recovery map with minimal worst-case error.
It is given by regularization via the map
\be
\label{RegProg}
\Delta^\sigma:
y \in \bR^m
\mapsto
f_{y}^\sigma := \Big[
\underset{f \in F}{\argmin \;} 
(1-\sigma) \|T(f)\|^2 + \sigma \|U(y-\La f)\|^2
\Big] \in F,
\ee
which is incidentally a linear map that can be expressed as $\Delta^\sigma =  [(1-\sigma) T^* T + \sigma \La^* U^* U \La]^{-1} \sigma \La^* U^* U$.
\cite{melkman1979optimal}
were the first to show that a globally optimal recovery map occurs as $\Delta^\sigma$ for {\em some} parameter $\sigma \in [0,1]$.
Much later, in the case $\Gamma  = \Id$, $T = \eps^{-1} P_{V^\perp}$, and $U = \eta^{-1} \Id$,
\cite{foucart2023optimal} exhibited a way to compute such an optimal regularization parameter by solving a semidefinite program.
This was extended to arbitrary (linear) $\Gamma$ and $T$ in \citep{foucart2023optimalSampTA}
and further simplified in \citep{foucart2024radius}, 
which also covered the case of an arbitrary (linear)~$U$.
In short, a globally optimal recovery map is given by $\Gamma \circ \Delta^{\sigma^\sharp}$,
where the optimal parameter is $\sigma^\sharp = b^\sharp/(a^\sharp + b^\sharp )$ with
\be
\label{NoisyMGWCEGamma}
(a^\sharp,b^\sharp)
= \underset{a,b \ge 0}{\argmin \;} a  + b 
\quad \mbox{s.to } \;
a \|Tf\|^2 + b \|U \La f\|^2 \ge \| \Gamma f\|^2
\; \; \forall f \in F.
\ee
When the Hilbert space $F$ is of finite dimension, the above can be solved as a semidefinite program,
since the constraint reformulates as
$$
a \, T^* T + b \, \La^* U^* U \La \succeq \Gamma^* \Gamma.
$$
Note that this semidefinite program holds the value of the minimal global worst-case error, 
namely
\be
\label{MinGWCEasSDP}
\min_\Delta {\rm gwce}(\Delta)^2 
= \min_{a,b \ge 0} \; a + b 
\quad \mbox{s.to } 
a \, T^* T + b \, \La^* U^* U \La \succeq \Gamma^* \Gamma.
\ee

In the particular situation
where $T = \eps^{-1 }P_{V^\perp}$, $U = \eta^{-1} \Id$,
and  $\La \La^* = \Id$---which occurs in our tensorial situation---any $\sigma \in [0,1]$ actually leads to optimality when $\Gamma = \Id$, i.e., for full recovery,
the global worst-case error of $\Delta^\sigma$ is independent of~$\sigma$, as shown in \citep{foucart2023optimal},
so there is no need to solve any semidefinite program.
Although this fact does not generalize to arbitrary quantities of interest,
the next subsection shall uncover other natural $\Gamma$ displaying this~peculiarity.

Another bonus of the particular situation
where $T = \eps^{-1 }P_{V^\perp}$, $U = \eta^{-1} \Id$,
and  $\La \La^* = \Id$
is a nice expression of the regularizer $f_{y}^{\sigma}$.
Namely, with the notation already used in  \eqref{FlaGram}, we have 
\be
\label{FlaGramInac}
f_{y}^{\sigma}
= \tau \sum_{i=1}^m a_i u_i + \sum_{j=1}^k b_j v_j,
\ee
where $\tau \in [0,1]$ depends on $\sigma \in [0,1]$ according to the correspondence
\be
\label{CorrespSigTau}
\tau = \f{\sigma \eps^2 }{(1-\sigma)\eta^2 + \sigma \eps^2}
\quad \longleftrightarrow \quad
\sigma = \f{\tau \eta^2}{(1-\tau) \eps^2 + \tau \eta^2}.
\ee
As a matter of fact,
it was established in \citep[Proposition 3 and Appendix]{foucart2023optimal} that the right-hand side of \eqref{FlaGramInac}
equals the minimizer of $(1-\tau) \|P_{V^\perp} f\|^2 + \tau \|y - \La f\|^2$,
which coincides with $f_{y}^{\sigma}$ under the correspondence \eqref{CorrespSigTau} when $T = \eps^{-1 }P_{V^\perp}$ and $U = \eta^{-1} \Id$.

\paragraph{Minimized local worst-case error.}
 
The {\em local} setting is the most delicate one,
only solved for the particular situation
$\Gamma = \Id$, $T = \eps^{-1} P_{V^\perp}$, $U = \eta^{-1} \Id$, and $\La \La^* = \Id$ in \citep{foucart2023optimal}---fortuitously, these restrictions hold in our tensor completion problem.
Fixing $y$, \cite{beck2007regularization} exhibited,
when\footnote{Technically speaking, the restriction $U = \eta^{-1} \Id$ was also in place,
but the case of an arbitrary $U$ can always be reduced to the case $U = \Id$ in the local setting,
since the set $\cS_U = \{f \in F: \|Tf\| \le 1, \|U(y - \La f)\| \le 1 \}$ of model- and data-consistent objects can be interpreted as $\cS_{\Id}$ by replacing $y$ by $Uy$ and $\La$ by $U \La$.} 
$\Gamma = \Id$,
a `candidate' Chebyshev center by uncovering an upper bound for ${\rm lwce}_y(g)^2$ and then minimizing over all $g$.
This is in the spirit of~\eqref{LWCE^2},
except that the S-lemma is replaced by the S-procedure,
which is not guaranteed to be exact,
so the relaxation only provides an upper bound.
Such a process was somewhat duplicated 
in Theorem~2 of \citep{foucart2024s} for an arbitrary (linear) $\Gamma$, leading to
\begin{align*}
\min_g {\rm lwce}_y(g)^2 \le 
\min_{\substack{a,b \ge 0\\ t \in \bR}} \;
a  + b   - b \|U y\|^2 + t
& & & \mbox{s.to } \,
a T^* T + b \La^* U^* U \La \succeq \Gamma^* \Gamma\\
& & & \mbox{and } \,
\bbmx
a T^* T + b \La^* U^* U \La & \vline & b \La^* U^* U y\\
\hline 
b (\La^* U^* U y)^*  & \vline & t
\ebmx \succeq 0,
\end{align*}
while the minimizers $a^\sharp,b^\sharp \ge 0$
yield the candidate Chebsyshev center $g^{\rm cand} = \Gamma(f_{y}^{\sigma^\sharp})$,
with $f_{y}^{\sigma^\sharp}$ still denoting the solution to the regularization problem from \eqref{RegProg} with parameter $\sigma^\sharp = b^\sharp/ (a^\sharp + b^\sharp)$.
The article \citep{foucart2024s} also indicated that
this candidate Chebyshev center agrees with the true Chebyshev center determined in \citep{foucart2023optimal} 
when $\Gamma = \Id$, $T = \eps^{-1} P_{V}^\perp$, $U = \eta^{-1} \Id$, and $\La \La^* = \Id$.
For arbitrary $\Gamma$, $T$, $U$, and $\La$,
\cite{foucart2025worst} later showed (see Theorem~5)
that the candidate Chebyshev center coincides with the true Chebyshev center if and only if the orthogonality condition 
$\langle T f_y^{\sigma^\sharp}, T h_{\sigma^\sharp} \rangle = 0$
hold,
where $h_\sigma$ generically denotes a leading eigenvector of $[(1-\sigma) T^* T + \sigma \La^* U^* U \La]^{-1/2} \Gamma^* \Gamma [(1-\sigma) T^* T + \sigma \La^* U^* U \La]^{-1/2}$.
Incidentally, these conditions were precisely the ones used as sufficient conditions 
in \citep{foucart2023optimal} 
(see the proof of Theorem~8 there) 
to determine the true Chebyshev center
in the specific situation 
where $\Gamma = \Id$, $T = \eps^{-1} P_{V^\perp}$, $U = \eta^{-1} \Id$, and $\La \La^* = \Id$.

All in all,
in this specific situation,
the true Chebyshev center can be computed as a regularizer~$f_y^{\sigma^\sharp}$ for some parameter $\sigma^\sharp \in [0,1]$
than can be determined by solving a semidefinite program.
But resorting to semidefinite programming is unnecessary:
working with the alternative parametrization \eqref{CorrespSigTau},
it was established in \citep{foucart2023optimal} that the optimal $\tau^\sharp \in [0,1]$ is the unique solution between $1/2$ and $\eps/(\eps+\eta)$ to the equation
\be
\label{EqMinEval=}
{\rm eig}_{\min}((1-\tau) P_{V^\perp} + \tau \La^* \La) = 
\f{(1-\tau)^2 \eps^2 - \tau^2 \eta^2}{(1-\tau) \eps^2 - \tau \eta^2 + (1-\tau) \tau (1-2\tau) \delta^2},
\ee
where $\delta$ is precomputed as the common (when $\La \La^* = \Id$) value of  $\min\{ \|P_{V^\perp} f\|: \La f = y\}$  and of $\min\{ \|\La f - y\|: f \in V \}$.
This equation can be---and was---solved efficiently via Newton method, 
since the derivatives of both sides are available.
It was indeed uncovered in the appendix of \citep{foucart2023optimal} that 
$e(\tau) := {\rm eig}_{\min}((1-\tau) P_{V^\perp} + \tau \La^* \La)$ satisfies
$$
\f{d e(\tau)}{d\tau}
= \f{1-2\tau}{\tau (1-\tau)} \f{e(\tau)(1-e(\tau))}{1-2 e(\tau)}.
$$
As such,
a possibly expensive eigencomputation is, 
at first sight,
needed at each Newton iteration.
This downside shall be removed in the next subsection.

After this step,
the value of the (squared) Chebyshev center becomes available, too, in the particular situation.
It was indeed revealed in 
\citep[last remark before Section 4]{foucart2023optimal}, 
and more generally in 
\citep{foucart2025worst}, 
that
\begin{align}
\label{MinLWCE=fla}
\min_g {\rm lwce}_y(g)^2
& = \f{(1-\tau^\sharp) \eps^2 + \tau^\sharp \eta^2 - (1-\tau^\sharp) \tau^\sharp \delta^2} {e(\tau^\sharp)}\\
\nonumber
& = \f{(1-\tau^\sharp) (\eps^2  -\|P_{V^\perp} f^\sharp \|^2 ) + \tau^\sharp (\eta^2 - \|y - \La f^\sharp\|^2 )} {e(\tau^\sharp)},
\end{align}
where the second expression is closer in spirit to \eqref{CheRadEigMax}.

\subsection{Two novel results}

The results presented in this section are both proved under the assumptions that $T$ is a multiple of an orthogonal projector,
that $U$ is a multiple of the identity,
and that $\La$ satisfies $\La \La^* = \Id$,
which means that $\La^* \La$ is an orthogonal projector from $F$ onto $\ran(\La^*)$.
These assumptions
are met in our tensor completion problem,
where $T = \eps^{-1} P_{V_s^\perp}$ and
$U = \eta^{-1} \Id$.
Moreover, it is readily checked that $\La \La^* = \Id$
because the observation map
$\La: \bR^{n_1 \times \cdots \times n_d} \to \bR^{\cO}$ for which $\La(X) = X_{|\cO}$
has adjoint $\La^*: \bR^{\cO} \to \bR^{n_1 \times \cdots \times n_d}$ given by
$\La^*(y) = \sum_{\mathbf{i} \in \cO} y_{\mathbf{i}} E_{\mathbf{i}}$,
where $E_{\mathbf{i}} \in \bR^{n_1 \times \cdots \times n_d}$ is the canonical tensor equal to one on the $\mathbf{i}$th entry and to zero on all other entries.
As mentioned in the previous subsection,
under these assumptions,
every regularization map $\Delta^\sigma$, $\sigma \in [0,1]$,
is globally optimal for the recovery of $\Gamma = \Id$.
We provide below a much simpler proof than the original one from \citep{foucart2023optimal}
and, in the process,
we uncover a similar result for further quantities of~interest.

\bthm
\label{ThmNov1}
Given the model and uncertainty sets $\{ f \in F: \|P_{V^\perp} f \| \le \eps \}$ and $\{ e \in \ell_2^m: \|e\| \le \eta \}$,
if $\La \La^* = \Id$,
then, 
when ${\rm range}(\Gamma^*) \inc \ker(\La)$ 
or when $\Gamma = \Id$,
one has
$$
{\rm gwce}(\Gamma \circ \Delta^\sigma) 
= \inf_{\Delta} {\rm gwce}(\Delta)
\qquad \mbox{ for all } \sigma \in (0,1),
$$
i.e.,
all regularization maps \eqref{RegProg},
when composed with $\Gamma$, produce globally optimal recovery maps.
\ethm

\bpf
For now, let $\Gamma$ be an arbitrary (linear) quantity of interest.  
According to \eqref{GWCE-lin} and 
\eqref{MinGWCEasSDP},
with the change $c = a \eps^{-2}$ and $d = b \eta^{-2}$,
both ${\rm gwce}(\Gamma \circ \Delta^\sigma)^2$ 
and $\inf_{\Delta} {\rm gwce}(\Delta)^2$
are mininizers of $c \eps^2 + d \eta^2$ subject to, respectively,
\begin{align}
\label{Const1}
c \|P_{V^\perp}f\|^2 + d  \| e\|^2 \phantom{\La}
& \ge \|\Gamma(f - \Delta^\sigma(\La f +e)) \|^2
& & \forall f\in F, e \in \bR^m,\\
\label{Const2}
c \|P_{V^\perp}f\|^2 + d \|\La f\|^2 
& \ge \| \Gamma f\|^2
& & \forall f \in F.
\end{align}
Thus,
it is enough to show that the constraints \eqref{Const1} and \eqref{Const2} are equivalent.
Since it is clear that \eqref{Const1} implies \eqref{Const2} by taking $e = - \La f$,
it remains to show that \eqref{Const2} implies \eqref{Const1}.
So let us assume that \eqref{Const2} holds in the form
\be
\label{Const2bis}
c \|P_{V^\perp}f\|^2 + d \|P_{W^\perp} f\|^2 
 \ge \| \Gamma f\|^2
\qquad \forall f \in F,
\ee
where $P_{W^\perp} = \La^* \La$ is the orthogonal projector onto $\ran(\La^*) = W^\perp$,
$W := \ker(\La)$.
Let us now consider $f \in F$ and $e \in \bR^m$.
Introducing $g:= f - \Delta^\sigma(\La f +e)$, 
we need to prove that 
\be
\label{Const1bis}
c\|P_{V^\perp} f \|^2 + d \|e\|^2 \overset{?}{\ge} \|\Gamma g \|^2.
\ee
Note that, in view of $\|e\|^2 = \| \La^* e\|^2$, the left-hand side is 
$$
c\|P_{V^\perp} f \|^2 + d \|e\|^2 
= \alpha + \beta,
\qquad \mbox{where }
\begin{cases}
\alpha = c \|P_{\phantom{(}V+W \phantom{)^\perp}} P_{V^\perp} f\|^2 + d \|P_{\phantom{(}V+W \phantom{)^\perp}}\La^* e\|^2,\\
\beta = c \|P_{(V+W)^\perp} P_{V^\perp} f\|^2 + d \|P_{(V+W)^\perp}\La^* e\|^2.
\end{cases}
$$
Because we can express the regularization map as\footnote{The existence of the inverse owes to the running assumption \eqref{NecCond} that $V\cap W = \{0\}$, which, in infinite dimensions, must be strengthened to the existence of some $\delta>0$ such that $\max \{ \|P_{V^\perp} f\|, \| P_{W^\perp} f\| \} \ge \delta \|f\|$ for all $f \in F$.}
$\Delta^\sigma = [(1-\sigma) P_{V^\perp} + \sigma P_{W^\perp}]^{-1} \sigma \La^*$,
we have
\begin{align*}
g & = [(1-\sigma) P_{V^\perp} + \sigma P_{W^\perp}]^{-1} 
((1-\sigma) P_{V^\perp}f + \sigma P_{W^\perp}f - \sigma \La^* (\La f + e))\\
& = [(1-\sigma) P_{V^\perp} + \sigma P_{W^\perp}]^{-1} 
((1-\sigma) P_{V^\perp} f - \sigma \La^* e). 
\end{align*}
This implies that $(1-\sigma) P_{V^\perp} g + \sigma P_{W^\perp}g = (1-\sigma) P_{V^\perp} f - \sigma \La^* e$,
and hence
\be
\label{in(V+W)perp} 
(1-\sigma) (P_{V^\perp} g - P_{V^\perp} f) 
= - \sigma (\La^* e + P_{W^\perp} g)
\quad \mbox{ belongs to }
 V^\perp \cap W^\perp = (V+W)^\perp.
\ee
Applying $P_{V+W}$ to \eqref{in(V+W)perp} leads to
$P_{V+W} (P_{V^\perp} g - P_{V^\perp} f) = 0 $
and $P_{V+W}(\La^* e + P_{W^\perp} g) = 0$, 
so that 
$$
P_{V+W}  P_{V^\perp} f = P_{V+W} P_{V^\perp} g 
\qquad \mbox{and} \qquad
P_{V+W}\La^* e = - P_{V+W} P_{W^\perp} g.
$$
As a result, we obtain
$$
\alpha = c \|P_{V+W } P_{V^\perp} g\|^2 + d \|P_{V+W } P_{W^\perp} g\|^2
= c \| P_{V^\perp} P_{V+W }  g\|^2 + d \|P_{W^\perp} P_{V+W } g\|^2,
$$
where the last step used the fact that orthogonal projectors $P_{Z'}$ and $P_{Z^\perp}$ commute for subspaces $Z \inc Z'$
(by virtue of $P_{Z^\perp} = \Id - P_{Z}$
and of $P_{Z'} P_Z = P_Z P_{Z'} = P_Z$).
Therefore, as a consequence of the constraint \eqref{Const2} expressed as in \eqref{Const2bis}, we arrive at 
$$
\alpha \ge \| \Gamma P_{V+W}  g\|^2.
$$
This is enough to brush aside the case ${\rm range}(\Gamma^*) \inc \ker(\La)$,
as then ${\rm range}(\Gamma^*) \inc V + W$ yields $P_{V+W} \Gamma^* = \Gamma^*$,
or equivalently $\Gamma P_{V+W}= \Gamma$,
so that the previous inequality reads $\alpha \ge \|\Gamma g \|^2$
and the inequality $\alpha + \beta \ge \|\Gamma g \|^2$ 
requested in \eqref{Const1bis} immediately follows.

For the quantity of interest $\Gamma = \Id$,
in addition to $\alpha \ge \|P_{V+W}g\|^2$,
we shall also establish that $\beta \ge \|P_{(V+W)^\perp}g\|^2 $,
yielding $\alpha + \beta \ge \|g\|^2$ to validate \eqref{Const1bis}.
To this end, we apply $P_{(V + W)^\perp}$ to \eqref{in(V+W)perp},
which leads to, in view of $(V + W)^\perp \inc V^\perp$ and $(V + W)^\perp \inc W^\perp$ and after some rearrangement,
$$
P_{(V+W)^\perp} g
= 
(1-\sigma) P_{(V+W)^\perp} f - \sigma P_{(V+W)^\perp}\La^* e.
$$
By convexity, we deduce that 
$$
\|P_{(V+W)^\perp} g \|^2
\le (1-\sigma) \|P_{(V+W)^\perp} f\|^2 + \sigma \| P_{(V+W)^\perp}\La^* e\|^2.
$$
Notice that the constraint \eqref{Const2bis} with $\Gamma = \Id$ forces $c,d \ge 1$ by taking some nonzero $f \in W$ and $f \in V$.
Consequently,
we arrive at 
$$
\|P_{(V+W)^\perp} g \|^2
\le c \|P_{(V+W)^\perp} f\|^2 + d \| P_{(V+W)^\perp}\La^* e\|^2
= c \|P_{(V+W)^\perp} P_{V^\perp} f\|^2 + d \| P_{(V+W)^\perp}\La^* e\|^2 = \beta.
$$
As this was the last inequality to establish,
the proof is now complete.
\epf

In our tensor completion problem,
suitable quantities of interest include restrictions to unobserved entries,
i.e., cases for which $\Gamma(X) = X_{| \cU}$ with $\cU \inc \cO^c$,
as then $\ran(\Gamma^*) = \spa\{ E_{\mathbf{i}}, \mathbf{i}  \in \cU \}$ is clearly included in $\ker(\La)$.
In these cases, there is even a more direct argument exploiting~\eqref{FlaGramInac}.
When written as $\Delta^\sigma (y) = \tau \sum_{\mathbf{i} \in \cO} a_i E_{\mathbf{i}} + \sum_{j=1}^{k_s} b_j v_j$,
this identity yields $\Gamma ( \Delta^\sigma (y))  = \sum_{j=1}^{k_s} b_j v_j$ independently of $\sigma$,
and hence ${\rm gwce}(\Gamma \circ \Delta^{\sigma}) = \sup \{ \| \Gamma(f) - \Gamma (\Delta^{\sigma}(\La f + e)) \|: \|P_{V_s^\perp} f \| \le \eps, \, \|e\| \le \eta \}$ is independent of $\sigma$, too.

For the second novel result,
we transition from the global setting to the local setting
and provide a more explicit solution to the full recovery problem in the situation where $T=\eps^{-1} P_{V^\perp}$, $U=\eta^{-1} \Id$, and $\La \La^* = \Id$.

\begin{thm}
\label{ThmNov2}
Given the model and uncertainty sets $\{ f \in F: \|P_{V^\perp} f \| \le \eps \}$ and $\{ e \in \ell_2^m: \|e\| \le \eta \}$,
if $\La \La^* = \Id$,
then the locally optimal recovery map for $\Gamma = \Id$ is given by
$$
y \in \bR^m
\mapsto
\Big[ 
\underset{f \in F}{\argmin \;} 
(1-\tau^\sharp) \|P_{V^\perp}f\|^2 + \tau^\sharp \|y-\La f\|^2
\Big] \in F,
$$
where the regularization parameter $\tau^\sharp$,
depending on $y$, 
is the unique solution between $1/2$ and $\eps/(\eps+\eta)$ to the sextic equation 
\be
\label{Sextic}
c \, \big( (1-\tau)\eps^2 - \tau \eta^2 + (1-\tau)\tau(1-2\tau)\delta^2 \big)^2
- (\eps^2 - \eta^2 + (1-2\tau) \delta^2) ((1-\tau)^2 \eps^2 - \tau^2 \eta^2) = 0,
\ee
with $c := {\rm eig}_{\min}(P_{V^\perp} + \La^* \La)(2-{\rm eig}_{\min}(P_{V^\perp} + \La^* \La)$ and $\delta := \min_{\La f = y} \|P_{V^\perp} f\|  = \min_{f \in V }  \|\La f - y\|$.
\end{thm}

\bpf
As recalled in the previous subsection,
the locally optimal recovery map is 
$y \in \bR^m \mapsto \Delta^{\sigma^\sharp} y$,
where the regularization parameter $\sigma^\sharp$
corresponds via \eqref{CorrespSigTau}
to the unique solution~$\tau^\sharp$ between $1/2$ and $\eps/(\eps+\eta)$ 
to
$$
e(\tau) = r(\tau),
\qquad \mbox{where} \quad
r(\tau) := \f{(1-\tau)^2 \eps^2 - \tau^2 \eta^2}{(1-\tau) \eps^2 - \tau \eta^2 + (1-\tau) \tau (1-2\tau) \delta^2}.
$$
Here $e(\tau) := {\rm eig}_{\min}((1-\tau) P_{V^\perp} + \tau \La^* \La)$ satisfies the differential equation
\be
\label{DiffEq}
\f{d e(\tau)}{d\tau}
= \f{1-2\tau}{\tau (1-\tau)} \f{e(\tau)(1-e(\tau))}{1-2 e(\tau)}.
\ee
The original article \citep{foucart2023optimal} failed to notice, however,
that this differential equation could be solved somewhat explicitly.
To this end,
consider the function of $\tau \in [0,1]$ defined by
$$
C(\tau) = \f{e(\tau)(1-e(\tau))}{\tau(1-\tau)}
= \f{e(\tau)-e(\tau)^2}{\tau - \tau^2},
$$
whose derivative, thanks to \eqref{DiffEq}, becomes 
$$
\f{dC(\tau)}{d\tau}
= 
\f{1-2e(\tau)}{\tau(1-\tau)} \f{de(\tau)}{d\tau} 
- \f{e(\tau)(1-e(\tau))}{(\tau(1-\tau))^2} (1-2\tau)
= 0.
$$
This implies that the function $C$ is constant,
equal to $c:=C(1/2)$, say, and we deduce that
\label{quadratic_tau}
$$
e(\tau)(1-e(\tau)) = c \, \tau(1-\tau)
\qquad \mbox{for all } \tau \in [0,1].
$$
On the one hand, since the optimal parameter $\tau^\sharp$,
lying between $1/2$ and $\eps / (\eps + \eta)$,
must satisfy $e(\tau) = r(\tau)$ and in turn 
$r(\tau)(1-r(\tau)) = c \, \tau(1-\tau)$,
we derive,
after expressing $1-r(\tau)$ as a ratio,
that $\tau^\sharp$ must be a solution to 
$$
\frac{\tau (1-\tau) (\eps^2 - \eta^2 + (1-2\tau) \delta^2) ((1-\tau)^2 \eps^2 - \tau^2 \eta^2)}{\big( (1-\tau)\eps^2 - \tau \eta^2 + (1-\tau)\tau(1-2\tau)\delta^2 \big)^2} = c \, \tau (1-\tau),
$$
which is readily equivalent to the announced sextic equation. 

On the other hand, let $\tau$ be a solution to this sextic equation lying between $1/2$ and $\eps / (\eps + \eta)$.
Reverting the above,
we obtain that $r(\tau)(1-r(\tau)) = c \, \tau(1-\tau)$,
and since $e(\tau)(1-e(\tau)) = c \, \tau(1-\tau)$,
we deduce that $e(\tau)(1-e(\tau)) =  r(\tau)(1-r(\tau))$.
This implies that
$e(\tau)=r(\tau)$---which is what we need for $\tau = \tau^\sharp$ to follow---or $e(\tau)=1-r(\tau)$---which we are aiming to invalidate.
To this end, it~is useful to recall from \citep[Proof of Theorem~8]{foucart2023optimal}
that the denominator $d(\tau)$ of $1-r(\tau)$ does not vanish between $1/2$ and $\eps/(\eps+\eta)$,
that $\delta \le \eps + \eta$,
and that $e(\tau) \in [0,1/2]$.
By contradiction,
assume now that $e(\tau) = 1-r(\tau) \not= r(\tau)$, so that $1-r(\tau) \in [0,1/2)$,
and consider e.g. the case $\eps \ge \eta$.
Since then $1/2 \le \tau \le \eps/(\eps+\eta)$,
the numerator $n(\tau)$ of $1-r(\tau)$ satisfies
$$
\f{n(\tau)}{\tau(1-\tau)} 
= \eps^2 - \eta^2 + (1-2\tau) \delta^2
\ge \eps^2 - \eta^2 + (1-2\tau) (\eps + \eta)^2
\ge \eps^2 - \eta^2 + \Big(1-\f{2 \eps}{\eps + \eta} \Big) (\eps + \eta)^2 = 0,
$$
i.e., its is nonnegative, 
hence the denominator $d(\tau)$ must be positive.
Therefore, the inequality $1 - r(\tau) < 1/2$ yields $2 n(\tau) - d(\tau) < 0$.
However, after a simplification step, we have
\begin{align*}
2 n(\tau) - d(\tau)
& = (2\tau-1) [(1-\tau) \eps^2 + \tau \eta^2 - \tau(1-\tau) \delta^2]\\
& \ge (2\tau-1) [(1-\tau) \eps^2 + \tau \eta^2 - \tau(1-\tau) (\eps^2 + \eta^2 + 2 \eps \eta)]\\
& = (2\tau-1) [(1-\tau)^2 \eps^2 + \tau^2 \eta^2 - 2 \tau(1-\tau)  \eps \eta]
= (2\tau - 1) [(1-\tau) \eps - \tau \eta]^2\\
& \ge 0.
\end{align*}
This provides the contradiction desired to finish the proof.
\epf

\subsection{Tensor computations}
\label{SubSecTensorCompInac}

As in the accurate scenario,
there are satisfying solutions to the abstract Optimal Recovery problem in the inaccurate scenario,
but further challenges occur for their computational implementations due to the high dimensionality of the tensor completion problem.
We discuss our resolutions below.

\paragraph{Optimal recovery maps.}
In both the global and local settings,
the optimal recovery maps are given by regularization.
When $T = \eps^{-1} P_{V^\perp}$ and $U = \eta^{-1} \Id$,
they take the form
$$
\Delta_\tau(y)
= f_{y,\tau}
= \underset{f \in F}{\argmin \;} 
(1-\tau) \|P_{V^\perp} f\|^2 + \tau \|y-\La f\|^2.
$$
Beware of the difference with $f^{\sigma}_y$ defined in \eqref{RegProg}---in fact, we have $f_{y,\tau} = f_y^{\sigma}$ under the correspondence~\eqref{CorrespSigTau}. 
Either way, an unconstrained least-squares problem needs to be solved.
Using $\tau$, the normal equations read
$A_\tau f_{y,\tau} = \tau \La^* y$,
where 
$A_\tau := (1-\tau)  P_{V^\perp} + \tau \La^* \La$ is a positive semidefinite operator---actually, positive definite thanks to our standing assumption \eqref{NecCond}.
This framework is favorable for conjugate gradient,
especially since the action of $A_\tau$ on a tensor $X$ splits onto two computationally cheap, matrix-free actions, namely: 
$P_{V^\perp}X = (\Id-\Pi_{\le s})X$ is computed via \eqref{CheapAction} in $\cO(n^d \ln(n^d))$ operations
and $\La^* \La X$, obtained from $X$ by zeroing out the entries outside of $\cO$, is computed in $\cO(m)$ operations.
We do not spell out the conjugate gradient iterations adapted to our case,
but we mention that each iteration involves two applications of $A_\tau$
and that the solution $f_{y,\tau}$ is achieved in a finite number $t_{\max}$ of iterations. 
It is guaranteed that $t_{\max} \le N = n_1 \cdots n_d$,
but $t_{\max}$ was typically much smaller in our experiments which imposed a stopping criterion based on residual.
All in all,
the total operation cost for computing the regularizer $f_{y,\tau}$ via conjugate gradient is $\cO(t_{\max} n^d \ln(n^d))$.

We note that $f_{y,\tau}$,
as the minimizer of $\|P_{V^\perp} f\|^2 + (\tau /(1-\tau) ) \|y-\La f\|^2$,
intuitively tend to the minimizer $f_y$ of $\|P_{V^\perp} f\|^2$ subject to $\La f = y$,
see \citep[Appendix]{foucart2024radius} for a proper justification of this statement.
The point here is that the Chebyshev center in the accurate scenario, i.e. $f_y$,
can be approximately computed via conjugate gradient, too.
In the particular case where $T = \eps^{-1 }P_{V^\perp}$, $U = \eta^{-1} \Id$,
and  $\La \La^* = \Id$,
it even occurs that $f_y = f_{y,0}$ can be computed {\em exactly} via two conjugate gradient calls.
Indeed, recalling the affine dependence \eqref{FlaGramInac} of $f_{y,\tau}$,
written as $f_{y,\tau} = f_{y,0} + \tau (f_{y,1} - f_{y,0})$,
we can deduce $f_y = f_{y,0}$ from $f_{y,\tau}$ computed at two values of $\tau \in (0,1)$,
e.g. $f_y = 2 f_{y,1/3} - f_{y,2/3}$.

To finally compute the optimal recovery maps,
now that we are able to compute each regularizer~$f_{y,\tau}$,
we need to determine the optimal parameter $\tau^\sharp$.
In the {\em global} setting,
this is generally done by solving the semidefinite program~\eqref{MinGWCEasSDP}, 
yet solving it naively will be inefficient due to the high dimensionality of our tensor completion problem.
Fortuitously, the assumptions $T = \eps^{-1} P_{V^\perp}$, $U = \eta^{-1} \Id$, and $\La \La^* = \Id$ hold for this problem, 
so any $\tau \in [0,1]$ leads to a globally optimal recovery map when $\Gamma(X) = X_{| \cU}$, $\cU \inc \cO^c$, and when $\Gamma(X)=X$,
as shown in Theorem~\ref{ThmNov1}.
In the {\em local} setting,
again thanks to the assumptions $T = \eps^{-1} P_{V^\perp}$, $U = \eta^{-1} \Id$, and $\La \La^* = \Id$ being met,
Theorem~\ref{ThmNov2} establishes that the optimal parameter is obtained by solving the sextic equation~\eqref{Sextic}. 
To do so,
we first need to compute $\delta$ and $c$.
The quantity $\delta = \min\{ \|P_{V^\perp} f\|: \La f = y\} = \min\{ \|\La f - y\|: f \in V \}$ is available after the computation of any $f_{y,\tau}$,
according to $\delta = \|P_{V^\perp} f_{y,\tau}\|/\tau = \|\La f_{y,\tau} - y \|/(1-\tau)$, see \citep[Proposition~3]{foucart2023optimal}.
The quantity $c = {\rm eig}_{\min}(P_{V^\perp} + \La^* \La)(2-{\rm eig}_{\min}(P_{V^\perp}+ \La^* \La)$ shall by determined by performing a power method only once.
For comparison,
before Theorem~\ref{ThmNov2},
the determination of the optimal parameter $\tau^\sharp$ necessitated one eigenanalysis (power method or else) at each Newton iteration.
To compute $c$, 
it is enough to compute 
$$
{\rm eig}_{\min}(P_{V^\perp} + \La^* \La)
= 2 - {\rm eig}_{\max}(B),
\qquad \mbox{where} \quad
B:= 2 \Id - P_{V^\perp} - \La^* \La \succeq 0.
$$
The standard power method to determine the largest absolute eigenvalue of $B$,
also applicable to matrices that are not positive semidefinite,
would produce the `powers' $x^{(t)} = B^t x^{(0)}$ of an arbitrary initial vector $x^{(0)}$
and approximate ${\rm eig}_{\max}(B)$ by $|\langle A x^{(t)}, x^{(t)} \rangle|/ \|x^{(t)}\|^2$.
We modify this strategy slightly by averaging over several random initial vectors
in order to provide a PAC (short for probably approximately correct) guarantee of accuracy with a number of iterations known beforehand. 
We emphasize that each iteration $x^{(t)} = B x^{(t-1)}$
still exploits the fact that the action of $B$ again splits as two computationally cheap, matrix-free actions.
The result we rely on is stated and proved below, although it is likely to have been expressed in a resembling form elsewhere before.

\bprop
\label{PropPower}
Let $B \in \bR^{N\times N}$ be positive semidefinite matrix with largest eigenvalue $\la_{\max}$.
Given integers $k,t \ge 0$,
consider independent standard gaussian vectors $g^{(1)},\ldots,g^{(k)} \in \bR^N$
and form the computable randomized estimator
$$
\wh{\la}_{k,t} :=
\bigg[
\frac{5}{k}
\sum_{j=1}^k \big\langle B^t g^{(j)}, g^{(j)} \big\rangle
\bigg]^{1/t}.
$$
Given two small parameters $\iota,\theta>0$,
as soon as $\displaystyle{k \ge 25 \ln \Big(  \dfrac{2}{\theta} \Big)}$
and
$\displaystyle{t \ge \dfrac{\ln(9N)}{\ln(1+\iota)}}$,
it holds with probability at least $1-\theta$ that
$$
\lambda_{\max}
\le
\wh{\la}_{k,t}
\le
(1+\iota)\lambda_{\max}.
$$
\eprop

\bpf
Write the eigendecomposition of $B$ as 
$$
B = \sum_{i=1}^N \la_i v_i v_i^*,
\quad \mbox{with } \la_{\max} = \la_1 \ge \cdots \ge \la_N \ge 0
\; \mbox{ and } \;
(v_1, \ldots,v_N) \mbox{ an orthonormal basis of } \bR^N. 
$$
For each standard gaussian vector $g^{(j)}$,
we have $g^{(j)} = \sum_{i=1}^N g^{(j)}_i v_i$,
where all the $g^{(j)}_i = \langle g^{(j)}, v_i \rangle$ are independent standard gaussian variables.
The random variable $\xi_{j,t}:=\langle B^t g^{(j)}, g^{(j)} \rangle = \sum_{i=1}^N \la_i^t (g^{(j)}_i)^2$
therefore has expectation 
$$
\bE_t: = \bE \big[ \langle B^t g^{(j)}, g^{(j)} \rangle \big]
= \sum_{i=1}^N \la_i^t \;
\begin{cases}
\le & N \la_{\max}^t,\\
\ge & \la_{\max}^t. \phantom{N}
\end{cases}
$$
Based on the moment generating function of the square of a standard gaussian variable (see e.g. \citep[(B.5)]{foucart2022mathematical}),
the moment generating function of $\xi_{j,t}$ at $u>0$ is
\begin{align*}
\bE[\exp(u \xi_{j,t})]
& = \prod_{i=1}^N \bE[\exp(u \la_i^t (g^{(j)}_i)^2)]
= \prod_{i=1}^N \frac{1}{\sqrt{1-2u\la_i^t}}\\
& \le \prod_{i=1}^N \exp( u \la_i^t + 4 u^2 \la_i^{2t} )
= \exp(u \bE_t + 4 u^2 \bE_{2t})\\
& \le \exp(u \bE_t + 4 u^2 \bE_{t}^2).
\end{align*}
As a result, for any $\omega \in (0,1)$
and making the choice $u = \omega/(8 \bE_t)$,
an application of Markov inequality yields
\begin{align*}
\bP \bigg[ \f{1}{k} \sum_{j=1}^k \big\langle B^t g^{(j)}, g^{(j)} \big\rangle & > (1+\omega) \bE_t \bigg]
= \bP \bigg[ \exp\Big( u \sum_{j=1}^k \xi_{j,t} \Big) 
> \exp(u k (1+\omega) \bE_t)
\bigg]\\
& \le \f{\bE \big[ \exp\big( u \sum_{j=1}^k \xi_{j,t} \big) \big]}{\exp(u k (1+\omega) \bE_t)}
= \prod_{j=1}^k  \f{\bE[\exp(u \xi_{j,t})]}{\exp(u(1+\omega) \bE_t)}
\le \prod_{j=1}^k  \f{\exp(u \bE_t + 4 u^2 \bE_{t}^2)}{\exp(u(1+\omega) \bE_t)}\\
& = \prod_{j=1}^k \exp(-u \omega \bE_t + 4 u^2 \bE_t^2)
= \prod_{j=1}^k \exp \Big( - \f{\omega^2}{16} \Big)
= \exp \Big( - \f{\omega^2 k}{16} \Big).
 \end{align*}
In a similar fashion and with details left to the reader,
we also derive
$$
\bP \bigg[ \f{1}{k} \sum_{j=1}^k  \big\langle B^t g^{(j)}, g^{(j)} \big\rangle < (1-\omega) \bE_t \bigg]
\le \exp \Big( - \f{\omega^2 k}{16} \Big).
$$
This means that, with failure probability at most $2 \exp(-\omega^2 k/16)$,
we have
$$
\f{1}{k} \sum_{j=1}^k  \big\langle B^t g^{(j)}, g^{(j)} \big\rangle \;
\begin{cases}
\le & (1+\omega) \bE_t \le (1+\omega)N \la_{\max}^t,\\
\ge & (1-\omega) \bE_t \ge (1-\omega) \la_{\max}^t, \phantom{N}
\end{cases}
$$
which implies that
$$
\la_{\max}
\le \bigg[ \f{1}{(1-\omega)k} \sum_{j=1}^k  \big\langle B^t g^{(j)}, g^{(j)} \big\rangle \bigg]^{1/t}
\le \bigg( \f{1+\omega}{1-\omega} N \bigg)^{1/t} \la_{\max}.
$$
By choosing $\omega = 4/5$,
we conclude that $\la_{\max} \le \wh{\la}_{k,t} \le (9N)^{1/t} \la_{\max}$ with failure probability at most $2 \exp(-k/25)$. 
We arrive at $\la_{\max} \le \wh{\la}_{k,t} \le (1+\iota) \la_{\max}$ with failure probability at most $\theta$ when $k \ge 25 \ln(2/\theta)$ and $t \ge \ln(9N)/\ln(1+\iota)$,
as announced in Proposition~\ref{PropPower}.
\epf

\paragraph{Minimal worst-case errors.}
In the general situation,
computing the minimal {\em global}  worst-case error via the semidefinite program~\eqref{MinGWCEasSDP} seems to be the only reliable grip to put our hands on.
But this approach is not suitable in the high-dimensional setting of tensors,
unless we probe the inner workings of semidefinite solvers in search for saving opportunities---a task which we have not attempted.
In the particular situation where  $T = \eps^{-1} P_{V^\perp}$, $U = \eta^{-1} \Id$, $\La \La^* = \Id$, and $\Gamma = \Id$,
although we know that the minimal global  worst-case error is achieved as ${\rm gwce}(\Delta^\sigma)$ for any the regularization parameter $\sigma \in [0,1]$,
the computation of each of these still materializes, at first sight,
as the  semidefinite program~\eqref{GWCE-lin}.
Fortunately,
it is also known that the minimal global worst-case error agrees with the minimal local worst-case error at $y=0$ in the general situation,
and since the latter can be efficiently computed at any $y$ in the particular situation---see below---the minimal global worst-case error will ensue 
(see \citep[last remark of Section~4]{foucart2023optimal} for details).

Now, in the local setting
and with $T = \eps^{-1} P_{V^\perp}$, $U = \eta^{-1} \Id$, $\La \La^* = \Id$,
and $\Gamma = \Id$,
we have described ways to efficiently compute the regularizers $f_{y,\tau}$,
as well as the optimal regularization parameters $\tau^\sharp$
(involving the precomputation of the quantity $\delta$, which depends on $y$).
The minimal eigenvalue $e(\tau^\sharp) = {\rm eig}_{\min}((1-\tau^\sharp) P_{V^\perp} + \tau^\sharp \La^* \La )$ follows by \eqref{EqMinEval=}
and finally the minimal local worst-case error at any $y$ is deduced via \eqref{MinLWCE=fla}.

\section{Numerical Experiments}
\label{sec:numerical-experiments}

In this section, we first verify
the accuracy and computational advantage of the FFT-based ANOVA projection.
Next, we compare global and local recoveries on synthetic data generated from
Sobol's product function. Finally, we evaluate our recovery methods on real-world data: the Energy Efficiency dataset serves as a moderate-scale example, 
while the Theoretical Multinary Perovskite Oxides dataset provides a large-scale test of computational scalability. 
{\sc matlab} files allowing one to 
retrieve our implementations and to
reproduce this section's experiments are available from the repository \url{https://github.com/Jingchun-Shao/OR_ANOVA}.

\subsection{Fast ANOVA via the fast Fourier transform}

We compare our FFT-based ANOVA construction with the traditional recursive construction on ten random tensors \(X\in\mathbb{R}^{n\times\cdots\times n}\). Here, with \(d\) denoting the order of the tensor, \(n\) the number of grid points in each dimension, and \(s\)  the maximum ANOVA interaction order, 
the experiments cover the range \(5\leq d\leq10\), \(2\leq n\leq4\), and \(1\leq s\leq4\).

\begin{figure}[ht]
\centering
\includegraphics[width=0.65\textwidth]{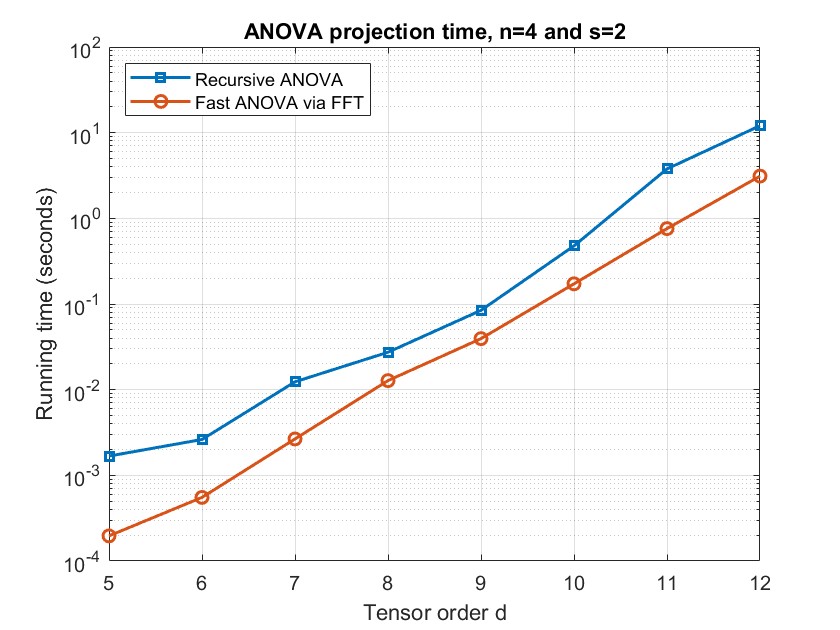}
\caption{Running times of the FFT-based and recursive ANOVA projections}
\label{fig:fast-anova-scaling}
\end{figure}

The relative error between the two implementations varies from \(3.09\times10^{-16}\) to \(2.19\times10^{-15}\), 
showing that they agree to machine precision.
We then compare their computational efficiencies. 
Figure~\ref{fig:fast-anova-scaling} compares running times as the order of the tensor increases from \(d=5\) to \(d=12\), with \(n=4\) and \(s=2\) being fixed. 
The FFT-based implementation is faster in every tested case. 
In addition, the dominant computational scaling of $n^d$ for both methods is reflected by the roughly linear trends seen with the logarithmic $y$-scale.

\subsection{Global and local recoveries for Sobol's function}
\label{SubSecSobol}

Motivated by Sobol's ANOVA-based sensitivity analysis \citep{sobol1993sensitivity} and Owen's effective-dimension analysis \citep{Owen2003}, we consider Sobol's product function defined on $[0,1]^d$ by
\[
x(t_1,\ldots,t_d)=\prod_{j=1}^{d}\frac{|4t_j-2|+a_j}{1+a_j}.
\]
We take \(d=8\), with \(a_1=a_2=0\) and \(a_3=\cdots=a_8=3\), 
so that the first two variables are more influential than the others six. 
Discretizing each variable at three points produces a tensor~$X$ with \(N=3^8=6561\) entries.
For \(s=1,2,3,4\), 
the relative ANOVA approximation errors $\|X - \Pi_{\le s} X\|/\|X\|$ are \(33.248\%\), \(10.218\%\), \(2.2759\%\), and \(0.3843\%\), respectively. 
We choose \(s=3\), yielding an approximation error of about \(2.3\%\) with model dimension \(k_s=577\).

The observations are taken at distinct entries sampled uniformly at random.
Their number is \(m=1968\),
making up \(30\%\) of the grid.
The vector of observations is gaussian with \(\|e\|=0.01\|\Lambda X\|\).
We choose the parameters $\eps$ and $\eta$
as \(\eps=1.1 \|X - \Pi_{\le s} X\|\)
and \(\eta=1.1\|e\|\).
In view of $\La \La^* = \Id$ in this situation,
the globally and locally optimal recovery procedures both return regularizers~$X_{\tau}$:
in the global setting, any $\tau \in [0,1]$ leads to the same global worst-case error,
so we choose $\tau = 0.5$;
in the local setting, the optimal parameter $\tau$ is found by solving the sextic equation. 
Note that finding it via the semidefinite formulation of \citet{beck2007regularization} was also attempted,
but the ambient dimension of $N = 6,561$
was too large for the computation to run given our resources.

Table~\ref{tab:sobol-recovery} reports the minimal local worst-case error computed as in \eqref{MinLWCE=fla},
as well as the minimal global worst-case error computed as a minimal local worst-case error with $y=0$,
both of them presented in normalized form,
i.e., after division by $\|X\|= 166.7529$.
This is possible in this synthetic example where the ground truth is available.
The local error is smaller than the global error, as should be.
Table~\ref{tab:sobol-recovery} also reports the pointwise, not worst-case, errors
$\|X-X_{\tau}\|$ (full grid)
and $\|X_{\cO^c} - (X_\tau)_{\cO^c} \|$ (unobserved grid), presented in relative forms,
i.e., divided by $\|X\|$ and $\|X_{\cO^c}\|$, respectively.
We see that the locally optimal procedure slightly reduces the full-grid relative error,
while the unobserved-grid relative error remains unchanged,
confirming the fact explained after Theorem~\ref{ThmNov1} that changing \(\tau\) only affects the observed entries in this situation.

\begin{table}[ht]
\centering
\begin{tabular}{c c c c c}
\hline
Method  & Normalized wce
& Full-grid relative error & Unobserved relative error \\
\hline
Global & \(11.999\%\) & \(2.5439\%\) & \(2.9301\%\) \\
Local  & \(11.174\%\) & \(2.5177\%\) & \(2.9301\%\) \\
\hline
\end{tabular}
\caption{Normalized worst-case errors and  pointwise relative errors for Sobol's function with \(s=3\).}
\label{tab:sobol-recovery}
\end{table}

\subsection{Global and local recoveries on the Energy Efficiency dataset}

The Energy Efficiency dataset of \citep{tsanas2012energy} comprises
\(768\) building designs and two responses: the heating load \(Y_1\) and
the cooling load \(Y_2\). 
The inputs \(X_1,\ldots,X_5\) describe the geometry,
\(X_6\) the orientation, 
and \(X_7,X_8\) the glazing area and its distribution.
In this dataset, each value of \(X_1\) actually determines the whole tuple \((X_1,\ldots,X_5)\), 
so that \(X_1\) alone indexes geometry.
We also treat each pair \((X_7,X_8)\) as one glazing configuration, giving \(16\)
glazing levels. 
In this way, the axes \(X_1\), \(X_6\), and the configurations \((X_7,X_8)\)
give rise to a \(12\times4\times16\) tensor with \(768\) entries. 
We randomly split these entries into a training set of \(538\) observed entries 
and a testing set of  \(230\) unobserved entries to be recovered.

For $s=0,1,2$, the relative ANOVA approximation errors are computed, respectively, as \(41.191\%\), \(4.0922\%\), \(1.2266\%\) for \(Y_1\), and \(36.064\%\), \(6.4641\%\),
\(4.8717\%\) for \(Y_2\). 
We retain $s=1$ and $s=2$
for our experiments, 
verifying in passing that the dimensions $k_1 = 30$ and $k_2 = 273$ do not exceed $m = 538$.
As for the model parameters $\eps$ and $\eta$,
we choose them in a heuristic way allowing us to deliberately avoid using the ground truth 
(even though it is available here).
We then proceed as in Subsection~\ref{SubSecSobol} to determine the normalized minimal worst-case errors and the testing (i.e., pointwise) relative errors.

Table~\ref{tab:energy-recovery} once again confirms that
the minimal local worst-case errors are smaller than their global counterparts
and 
that the local and global testing errors coincide because changing \(\tau\)
only affects the observed entries.
It also shows that
increasing \(s\) from \(1\)
to \(2\) reduces the testing error for \(Y_1\) from \(3.98\%\) to
\(2.48\%\)---both being excellent---but increases it from \(6.66\%\) to \(8.71\%\) for \(Y_2\).
Thus, on this split, second-order recovery predicts \(Y_1\) better,
while first-order recovery predicts \(Y_2\) better.

\begin{table}[ht]
\centering
\begin{tabular}{c c c c c c}
\hline
Target & \(s\) & Method
& Normalized wce & Testing relative error \\
\hline
\(Y_1\) & 1 & Global  & \(9.552\%\) & \(3.9782\%\) \\
\(Y_1\) & 1 & Local   & \(8.452\%\) & \(3.9782\%\) \\
\(Y_1\) & 2 & Global  & \(15.603\%\) & \(2.4844\%\) \\
\(Y_1\) & 2 & Local   & \(15.368\%\) & \(2.4844\%\) \\
\hline
\(Y_2\) & 1 & Global  & \(13.564\%\) & \(6.6573\%\) \\
\(Y_2\) & 1 & Local   & \(11.628\%\) & \(6.6573\%\) \\
\(Y_2\) & 2 & Global  & \(29.480\%\) & \(8.7145\%\) \\
\(Y_2\) & 2 & Local   & \(26.851\%\) & \(8.7145\%\) \\
\hline
\end{tabular}
\caption{Normalized worst-case errors and  testing relative errors for the Energy Efficiency dataset.
}
\label{tab:energy-recovery}
\end{table}

\subsection{Global recovery on the Multinary Perovskite Oxides dataset}

Discovering functional materials is an important but difficult task: 
the number of
candidate compositions grows rapidly, while evaluating each by
first-principles calculations is costly. 
\citet*{bare2023multinary} report formation enthalpies for \(66{,}516\)
multinary oxide compositions. Their dataset includes compositions with
two cations on both sites, \(A_{0.5}A'_{0.5}B_{0.5}B'_{0.5}O_3\), allowing \(A=A'\) or \(B=B'\)
when one site contains only a single cation. Motivated by the successful application of HDMR to molecular system  in \citep*{LiRabitz2001}, we adopt a low-order ANOVA model to capture the effects of individual cation choices and their interactions.

We use the four cation identities \(A_1,A_2,B_1,B_2\) as tensor coordinates,
giving a \(39\times38\times39\times38\) grid with \(N=2{,}196{,}324\)
entries. 
Of the recorded compositions, \(53{,}213\) are used for training
and \(13{,}303\) for testing. 
We recover the formation enthalpy 
via the globally optimal procedure by choosing
\(\tau=0.5\) and \(s=1,2\). 
The locally optimal procedure,
which is still too costly in these dimensions,
was not deployed.

Table~\ref{tab:multinary-oxides-recovery} 
reports the testing relative errors and
shows that including
second-order interactions reduces them from
\(4.33\%\) to \(3.38\%\). 
Both low-order models achieve test-set
relative errors below \(5\%\), highlighting the potential of our
method to support materials discovery.

\begin{table}[ht]
\centering
\begin{tabular}{c c c}
\hline
\(s\) & Model dimension \(k_s\) & Testing relative error \\
\hline
1 & 151  & \(4.330\%\) \\
2 & 8588 & \(3.379\%\) \\
\hline
\end{tabular}
\caption{Recovery of formation enthalpy for the Theoretical Multinary
Perovskite Oxides dataset.}
\label{tab:multinary-oxides-recovery}
\end{table}

\bibliography{refs}

@article{foucart2024radius,
  title={Radius of information for two intersected centered hyperellipsoids and implications in optimal recovery from inaccurate data},
  author={Foucart, S. and Liao, C.},
  journal={Journal of Complexity},
  volume={83},
  pages={101841},
  year={2024},
  publisher={Elsevier}
}

@article{foucart2025worst,
  title={Worst-Case Learning under a Multifidelity Model},
  author={Foucart, S. and Hengartner, N.},
  journal={SIAM/ASA Journal on Uncertainty Quantification},
  volume={13},
  number={1},
  pages={171--194},
  year={2025},
  publisher={SIAM}
}

@article{foucart2023optimal,
  title={Optimal recovery from inaccurate data in {H}ilbert spaces: regularize, but what of the parameter?},
  author={Foucart, S. and Liao, C.},
  journal={Constructive Approximation},
  volume={57},
  number={2},
  pages={489--520},
  year={2023},
  publisher={Springer}
}

@inproceedings{foucart2023optimalSampTA,
  title={On the optimal recovery of graph signals},
  author={Foucart, S. and Liao, C. and Veldt, N.},
  booktitle={2023 International Conference on Sampling Theory and Applications (SampTA)},
  pages={1--5},
  year={2023},
  organization={IEEE}
}

@article{beck2007regularization,
  title={Regularization in regression with bounded noise: A {C}hebyshev center approach},
  author={Beck, A. and Eldar, Y. C.},
  journal={SIAM Journal on Matrix Analysis and Applications},
  volume={29},
  number={2},
  pages={606--625},
  year={2007},
  publisher={SIAM}
}

@incollection{foucart2024s,
  title={S-procedure relaxation: A case of exactness involving {C}hebyshev centers},
  author={Foucart, S. and Liao, C.},
  booktitle={Explorations in the Mathematics of Data Science: The Inaugural Volume of the Center for Approximation and Mathematical Data Analytics},
  pages={1--18},
  year={2024},
  publisher={Springer}
}

@article{takemura1983tensor,
  title={Tensor analysis of {ANOVA} decomposition},
  author={Takemura, A.},
  journal={Journal of the American Statistical Association},
  volume={78},
  number={384},
  pages={894--900},
  year={1983},
  publisher={Taylor \& Francis}
}

@article{foucart2022learning,
  title={Learning from non-random data in {H}ilbert spaces: an optimal recovery perspective},
  author={Foucart, S. and Liao, C. and Shahrampour, S. and Wang, Y.},
  journal={Sampling Theory, Signal Processing, and Data Analysis},
  volume={20},
  number={1},
  pages={5},
  year={2022},
  publisher={Springer}
}

@article{binev2017data,
  title={Data assimilation in reduced modeling},
  author={Binev, P. and Cohen, A. and Dahmen, W. and DeVore, R. and Petrova, G. and Wojtaszczyk, P.},
  journal={SIAM/ASA Journal on Uncertainty Quantification},
  volume={5},
  number={1},
  pages={1--29},
  year={2017},
  publisher={SIAM}
}

@book{foucart2022mathematical,
  title={Mathematical Pictures at a Data Science Exhibition},
  author={Foucart, S.},
  year={2022},
  publisher={Cambridge University Press}
}

@article{melkman1979optimal,
  title={Optimal estimation of linear operators in {H}ilbert spaces from inaccurate data},
  author={Melkman, A. A. and Micchelli, C. A.},
  journal={SIAM Journal on Numerical Analysis},
  volume={16},
  number={1},
  pages={87--105},
  year={1979},
  publisher={SIAM}
}

@book{kowalski1995selected,
  title={Selected Topics in Approximation and Computation},
  author={Kowalski, M. A. and Sikorski, K. A. and Stenger, F.},
  year={1995},
  publisher={Oxford University Press}
}

@article{drautz_ace,
  title   = {Atomic cluster expansion for accurate and transferable interatomic potentials},
  author  = {Drautz, R.},
  journal = {Physical Review B},
  volume  = {99},
  number  = {1},
  pages   = {014104},
  year    = {2019},
  doi     = {10.1103/PhysRevB.99.014104}
}

@article{cohen2012capturing,
  title={Capturing ridge functions in high dimensions from point queries},
  author={Cohen, A. and Daubechies, I. and DeVore, R. and Kerkyacharian, G. and Picard, D.},
  journal={Constructive Approximation},
  volume={35},
  number={2},
  pages={225--243},
  year={2012},
  publisher={Springer}
}

@article{devore2011approximation,
  title={Approximation of functions of few variables in high dimensions},
  author={DeVore, R. and Petrova, G. and Wojtaszczyk, P.},
  journal={Constructive Approximation},
  volume={33},
  number={1},
  pages={125--143},
  year={2011},
  publisher={Springer}
}

@article{sloan1998quasi,
  title={When are quasi-{M}onte {C}arlo algorithms efficient for high dimensional integrals?},
  author={Sloan, I. H. and Wo{\'z}niakowski, H.},
  journal={Journal of Complexity},
  volume={14},
  number={1},
  pages={1--33},
  year={1998},
  publisher={Elsevier}
}

@article{steinlechner2016riemannian,
  title={Riemannian optimization for high-dimensional tensor completion},
  author={Steinlechner, M.},
  journal={SIAM Journal on Scientific Computing},
  volume={38},
  number={5},
  pages={S461--S484},
  year={2016},
  publisher={SIAM}
}

@article{JainOh2014,
  title={Provable tensor factorization with missing data},
  author={Jain, P. and Oh, S.},
  journal={Advances in Neural Information Processing Systems},
  volume={27},
  year={2014}
}

@article{LiuEtAl2013,
  title={Tensor completion for estimating missing values in visual data},
  author={Liu, J. and Musialski, P. and Wonka, P. and Ye, J.},
  journal={IEEE Transactions on Pattern Analysis and Machine Intelligence},
  volume={35},
  number={1},
  pages={208--220},
  year={2013},
  publisher={IEEE}
}

@article{KressnerSteinlechnerVandereycken2014,
  title={Low-rank tensor completion by {R}iemannian optimization},
  author={Kressner, D. and Steinlechner, M. and Vandereycken, B.},
  journal={BIT Numerical Mathematics},
  volume={54},
  number={2},
  pages={447--468},
  year={2014},
  publisher={Springer}
}

@article{Zhang2016Cross,
  title={Cross: Efficient low-rank tensor completion},
  author={Zhang, A. R.},
  journal={The Annals of Statistics},
  volume={47},
  number={2},
  pages={936--964},
  year={2019},
  publisher={Institute of Mathematical Statistics}
}

@article{HillarLim2013,
  title={Most tensor problems are {NP}-hard},
  author={Hillar, C. J. and Lim, L.-H.},
  journal={Journal of the ACM (JACM)},
  volume={60},
  number={6},
  pages={1--39},
  year={2013},
  publisher={ACM New York, NY, USA}
}

@article{Owen2003,
  title={The dimension distribution and quadrature test functions},
  author={Owen, A. B.},
  journal={Statistica Sinica},
  volume={13},
  number={1},
  pages={1--17},
  year={2003},
  publisher={Institute of Statistical Science, Academia Sinica}
}

@article{LiRabitz2001,
  title={General foundations of high-dimensional model representations},
  author={Li, G. and Wang, S.-W. and Rosenthal, C. and Rabitz, H.},
  journal={Journal of Mathematical Chemistry},
  volume={30},
  number={1},
  pages={1--30},
  year={2001},
  publisher={Springer}
}

@inproceedings{Griebel2006,
  title={Sparse grids and related approximation schemes for higher dimensional problems},
  author={Griebel, M.},
  booktitle={Foundations of Computational Mathematics, Santander 2005},
  pages={106--161},
  year={2006},
  publisher={Cambridge University Press}
}

@book{dudgeon1984multidimensional,
  title={Multidimensional Digital Signal Processing},
  author={Dudgeon, D. E. and Mersereau, R. M.},
  year={1984},
  publisher={Prentice-Hall},
  series={Prentice-Hall Signal Processing Series},
  address={Englewood Cliffs, NJ}
}

@article{sobol1993sensitivity,
  author  = {Sobol', I. M.},
  title   = {Sensitivity Estimates for Nonlinear Mathematical Models},
  journal = {Mathematical Modeling and Computational Experiment},
  volume  = {1},
  number  = {4},
  pages   = {407--414},
  year    = {1993}
}

@misc{tsanas2012energy,
  author={Tsanas, A. and Xifara, A.},
  title={{E}nergy {E}fficiency [{D}ataset]},
  year={2012},
  howpublished={UC Irvine Machine Learning Repository},
  url={https://archive.ics.uci.edu/ml/datasets/Energy+efficiency}
}

@article{bare2023multinary,
  author  = {Bare, Z. J. L. and Morelock, R. J. and Musgrave, C. B.},
  title   = {Dataset of theoretical multinary perovskite oxides},
  journal = {Scientific Data},
  volume  = {10},
  pages   = {244},
  year    = {2023}
}

\end{document}